\documentstyle[12pt,amscd]{amsart}
\newtheorem{Theorem}{Theorem}[section]
\newtheorem{Lemma}[Theorem]{Lemma}
\newtheorem{Corollary}[Theorem]{Corollary} 
\newtheorem{Proposition}[Theorem]{Proposition}

\def\Proj{\operatorname{Proj}}  \def\Spec{\operatorname{Spec}}
\def\height{\operatorname{ht}} \def\rdim{\operatorname{rdim}} \def\sdim{\operatorname{sdim}}\def\ara{\operatorname{ara}} 
\def\To{\longrightarrow} \def\sk{\par\smallskip} \def\ms{\medskip}

\begin{document}
\title{Positivity of mixed multiplicities}
\author{Ng\^o Vi\^et Trung} 
\address{Institute of Mathematics,  Box 631, B\`o H\^o, Hanoi, Vietnam}
\email{nvtrung@@hn.vnn.vn}
\thanks{The author is partially supported by the National Basic Research Program} 
\maketitle

\vspace{-0.4cm}
\centerline{\small\it Dedicated to the memory of my teacher Wolfgang Vogel\footnote{Wolfgang Vogel would become 60 years old in the year 2000}}\bigskip

\section*{Introduction}\sk

Let $R = \oplus_{(u,v) \in {\Bbb N}^2} R_{(u,v)}$ be a standard bigraded algebra over an artinian local ring $K = R_{(0,0)}$. Standard means $R$ is generated over $K$ by a finite number of elements of degree $(1,0)$ and $(0,1)$. Since the length $\ell(R_{(u,v)})$ of $R_{(u,v)}$ is finite, we can consider $\ell(R_{(u,v)})$ as a function in two variables $u$ and $v$. This function was first studied by van der Waerden [W] and Bhattacharya [B] who proved that there is a polynomial $P_R(u,v)$ of degree $\le \dim R-2$ such that $\ell(R_{(u,v)}) = P_R(u,v)$ for $u$ and $v$ large enough. Katz, Mandal and Verma [KMV] found out that the degree of $P_R(u,v)$ is equal to $\rdim R-2$, where $\rdim R$ is the relevant dimension of $R$ defined as follows. Let $R_{++}$ denote the ideal generated by the homogeneous elements of degree $(u,v)$ with $u \ge 1$, $v \ge 1$. Let $\Proj R$ be the set of all homogeneous prime ideals $\wp \not\supseteq R_{\tiny ++}$ of $R$.  Then 
$$\rdim R := \max\{\dim R/\wp|\ \wp \in \Proj R \}$$
if $\Proj R \neq \emptyset$ and $\rdim R$ can be any negative integer if $\Proj R = \emptyset$. If  we write
$$P_R(u,v) = \sum_{i+j \le \rdim R-2} {a_{ij}}{u \choose i}{v \choose j},$$
then $a_{ij}$ are non-negative integers for  $i+j = \rdim R - 2$. Let us denote $a_{ij}$ by $e_{ij}(R)$ for all $i,j \ge 0$ with $e_{ij}(R) = 0$ for $i+j > \rdim R -2$.\sk

We call $P_R(u,v)$ the {\it Hilbert polynomial} and the numbers $e_{ij}(R)$ with $i+j = \rdim R-2$ the {\it mixed multiplicities} of $R$. These notions seem to be of fundamental importance. But we would be surprised to learn how little is known on mixed multiplicities of an arbitrary standard bigraded algebra, especially on their positivity [KMV], [KR1], [HHRT]. \sk

The main motivation for the study of mixed multiplicities comes from the following situation. Let $(A, \frak m)$ be a local ring. Given an $\frak m$-primary ideal $I$ and an ideal $J$ of $A$, we can consider the function $\ell(I^vJ^u/I^{v+1}J^u)$ of the standard bigraded algebra 
$$R(I|J) := \oplus_{(u,v) \in {{\Bbb N}}^2}I^vJ^u/I^{v+1}J^u.$$
over the artinian local ring $A/I$. Let $r = \deg P_{R(I|J)}(u,v)$. Then we put
$$e_i(I|J) := e_{ir-i}(R(I|J))$$
and call it the $i$th mixed multiplicity of $I$ and $J$, $i = 0,\ldots,r$. The name goes back to Teissier's Carg\`ese paper [Te] on complex analytic hypersurfaces with isolated singularities which interpreted the Milnor numbers of general linear sections ($\mu^*$-invariant) as mixed multiplicities. \sk

If $J$ is an $\frak m$-primary ideal, Bhattacharya [B] proved that $r = d-1$, $d = \dim A$. Later Rees [R1] showed that $e_0(I|J)$ is the (Samuel's) multiplicity $e(I,A)$ of $A$ with respect to $I$. This result was generalized by Risler and Teissier [Te] who proved that $e_i(I|J)$ is the multiplicity of $A$ with respect to an ideal generated by $i$ sufficiently general elements of $J$ and $d-i$ sufficiently general elements of $I$. From this it follows that $e_i(I|J)$ is always positive.
Subsequently, Rees [R2] gave a more accessible approach by showing that $e_i(I|J)$ is the multiplicity of $A$ with respect to a joint reduction of ideals. Note that Rees, Risler and Teissier considered the function $\ell(A/I^vJ^u)$ which is a polynomial $Q(u,v)$ of degree $d$ for $u, v \gg 0$. If we write 
$$Q(u,v) = \sum_{i+j \le d} {b_{ij}} {u \choose i}{v \choose j},$$ 
then $b_{id-i} = e_i(I|J)$ for $i = 0,\ldots,d-1$ and $b_{d0} = e(J,A)$. However, if $J$ is not a $\frak m$-primary ideal, the function $\ell(A/I^vJ^u)$ has no meaning for $u > 0$. \sk

To extend Teissier's result to hypersurfaces with non-isolated singularities we need to consider the case $J$ is not a $\frak m$-primary ideal. But this case has remained mysterious. A characterization of $e_i(I|J)$ as the multiplicity of $A$ with respect to sufficiently general elements or to joint reductions of ideals as in the $\frak m$-primary case has not been known. Even the positivity of $e_i(I|J)$ is not well understood.
Katz and Verma [KV] proved that if $\height J > 0$, $e_0(I|J) = e(I,A)$, $e_i(I|J) > 0$ for $i < s(J)$ and $e_i(I|J) = 0$ for $i \ge s(J)$, where $s(J)$ denotes the analytic spread of $J$. However, there is a counter-example to the claim $e_i(I|J) > 0$ for $i < s(J)$ (see Section 3).
On the other hand, various results have made mixed multiplicities more interesting in the general case. For instance, Verma [V1], [V2] and Katz and Verma [KV] discovered that the multiplicities of the Rees algebra and the extended Rees algebra of $J$ can be expressed as sums of mixed multiplicities. P.~Roberts [Ro] described local Chern classes in terms of mixed multiplicities. Achilles and Manaresi [AM] interpreted the degree of St\"uckrad-Vogel intersection cycles in terms of mixed multiplicities. Moreover, mixed multiplicities of ideals have been extended to modules by Kirby and Rees [KR1], [KR2], Kleiman and Thorup [KT1], [KT2]. These works have led to the following\ms

\noindent{\bf Problem.} Which mixed multiplicity is positive and how to compute it effectively? \ms

This paper will solve this problem for both mixed multiplicities of bigraded algebras and of ideals. We shall see that a mixed multiplicity is positive if and only if a certain ring has maximal dimension and that the positive mixed multiplicities can be expressed as Samuel's multiplicities. Our main tool is the notion of filter-regular sequences in a standard bigraded algebra $R$ (see Section 1). This notion originated from the theory of Buchsbaum rings which have their roots in intersection theory [SV2].\sk

Now we are going to present the main results of this paper. For this we shall need the following notation. For any pair of ideals $\frak a, \frak b$ of a commutative ring $S$ let
$${\frak a}:{\frak b}^\infty := \{x \in S|\ \text{there is a positive integer $n$ such that $x{\frak b}^n \subseteq {\frak a}$}\}.$$ 

>From the formula $\deg P_R(u,v) = \rdim R-2$ we can easily deduce that $\deg P_R(u,v) = \dim R/0:R_{++}^\infty-2$. Using Northcott and Rees' theory of reductions of ideals [NR] we can also express the partial degrees $\deg_u P_R(u,v)$ and $\deg_vP_R(u,v)$ of $P_R(u,v)$ with respect to  $u$ and $v$ in terms of the dimension of certain factor rings of $R$. Let $R_{(1)}$ and $R_{(2)}$ denote the ideals of $R$ generated by the homogeneous elements of degree $(1,0)$ and $(0,1)$, respectively. Then
\begin{eqnarray*} \deg_uP_R(u,v) & = & \dim R/(0:R_{++}^\infty+R_{(2)})-1,\\ \deg_vP_R(u,v) & = & \dim R/(0:R_{++}^\infty+R_{(1)})-1. \end{eqnarray*}

For short let $r =  \dim R/0:R_{++}^\infty-2$. For the mixed multiplicities $e_{ij}(R)$, $i+j = r$, we obtain the following result which gives an effective criterion for their positivity and expresses the positive mixed multiplicities  as Samuel's multiplicities of graded algebras. \ms

\noindent{\bf Theorem \ref{nonzero}.} {\em Let $i, j$ be non-negative integers, $i+j = r$. Let $x_1,\ldots,x_i$ be a filter-regular sequence of homogeneous elements of degree $(1,0)$. Then $e_{ij}(R) > 0$ if and only if $$\dim R/((x_1,\ldots,x_i):R_{++}^\infty+R_{(1)}) = j+1.$$ In this case, if we choose homogeneous elements $y_1,\ldots,y_j$ of degree $(0,1)$ such that $x_1,\ldots,x_i,y_1,\ldots,y_j$ is a filter-regular sequence, then} 
$$e_{ij}(R) = e(R/(x_1,\ldots,x_i,y_1,\ldots,y_j):R_{++}^\infty).$$ \sk

We may replace the condition $x_1,\ldots,x_i,y_1,\ldots,y_j$ being a filter-regular sequence by the condition $x_1,\ldots,x_i,y_1,\ldots,y_j$ being sufficiently general elements. But the notion of sufficiently general elements is vague, whereas the filter-regular property can be tested effectively.\sk

Let $r_1 = \dim R/(0:R_{++}^\infty+R_{(2)})-1$ and $r_2 = \dim R/(0:R_{++}^\infty  + R_{(1)})-1$. Then $e_{ir-i}(R) = 0$ for $i > r_1$ or $i < r-r_2$. Katz, Mandal and Verma [KMV] showed that the mixed multiplicities $e_{ir-i}(R)$, $r-r_2 \le i \le r_1$, can be any sequence of non-negative integers with at least a positive entry. Surprisingly, using Grothendieck's Connectedness Theorem we can prove that these mixed multiplicities are positive under some mild conditions: \ms

\noindent{\bf Corollary \ref{rigid2}.} {\em Let $R$ be a domain or a Cohen-Macaulay ring. Then $e_{ir-i}(R) > 0$ for $r-r_2 \le i \le r_1$.}\ms

In the same vein we can study the positivity of mixed multiplicities of ideals. But we have to involve the standard bigraded algebra 
$$R(J|I) := \oplus_{(u,v) \in {{\Bbb N}}^2}I^vJ^u/I^vJ^{u+1}.$$

\noindent{\bf Theorem \ref{main}.} {\em Let $J$ be an arbitrary ideal of $A$ and $0 \le i < s(J)$. Let $a_1,\ldots,a_i$ be elements in $J$ such that their images in $J/IJ$ and $J/J^2$ form filter-regular sequences in $R(I|J)$ and $R(J|I)$, respectively. Then $e_i(I|J) > 0$ if and only if  $$\dim A/(a_1,\ldots,a_i):J^\infty = \dim A/0:J^\infty-i.$$ In this case, $e_i(I|J) = e(I,A/(a_1,\ldots,a_i):J^\infty)$.} \ms

This result provides an effective way to check the positivity of mixed multiplicities and to compute the positive mixed multiplicities of ideals. Since sufficiently general elements in $J$ satisfy the assumption of Theorem \ref{main}, we can easily derive from it Risler and Teissier's result on mixed multiplicities of $\frak m$-primary ideals. In particular, we obtain several new insights on the range of  positive mixed multiplicities.\ms

\noindent{\bf Corollary \ref{rigid3}.} {\em Let $\rho = \max \{i|\ e_i(I|J) > 0\}.$ Then\par
{\rm (i) } $\height J -1 \le \rho < s(J)$,\par
{\rm (ii) } $e_i(I|J) > 0$ for $0 \le i \le \rho$,\par
{\rm (iii)} $\max \{i|\ e_i(I'|J) > 0\} = \rho$ for any $\frak m$-primary ideal $I'$ of $A$.} \ms

The last two properties are especially interesting because they show that {\it the positive mixed multiplicities concentrate in a rigid range which does not depend on} $I$. We will give an  example with  $\rho  < s(J)-1$ and we will show that $\rho = s(J)-1$ is often the case. \ms

\noindent{\bf Corollary \ref{rigid4}.} {\em  Assume that $A/0:J^\infty$ satisfies the first chain condition. Then $e_i(I|J) > 0$ for $0 \le i \le s(J)-1$.}\ms

We will apply Theorem \ref{main} to study the case $I = \frak m$, where mixed multiplicities can be used to compute the Milnor number of analytic hypersurfaces, the multiplicity of the Rees algebras, and the degree of projective embeddings of rational $n$-folds obtained by blowing up projective spaces. As examples we show that the Milnor numbers of general linear sections of an analytic hypersurface are mixed multiplicities of the Jacobian ideal (this fact was proved by Teissier only for isolated singularies) and we will compute the mixed multiplicities $e_i({\frak m}|J)$ in the following cases: \par
(1) $J$ is the defining ideal of a set of points in ${\Bbb P}^2$ with $h$-vector of decreasing type (which arises as plane sections of curves in ${\Bbb P}^3$ [MR], [CO]),  \par
(2) $J$ is the defining prime ideal of a curve in ${\Bbb P}^3$ lying on the quadric $x_0x_3-x_1x_2$ (inspired by Huneke and Huckaba's work on symbolic powers of such ideals [HH]), \par
(3) $J$ is a homogeneous prime ideal of analytic deviation 1 which is generated by forms of the same degree in a polynomial ring. \sk

We would like to mention that there is in the references only one class of non $\frak m$-primary ideals $J$ for which $e_i({\frak m}|J)$ have been computed in terms of the usual multiplicities [RV]. They are ideals generated by certain quadratic sequence (a generalisation of Huneke's $d$-sequences). These ideals can be also handled by Theorem \ref{main}.\sk

Another interesting application of Theorem \ref{main} is the following description of the degree of the St\"uckrad-Vogel cycles by means of mixed multiplicities of ideals. These cycles were introduced in order to prove a refined Bezout's theorem [SV1], [Vo]. By a result of van Gastel [Ga] their rational components correspond to the distinguished varieties in Fulton's intersection theory.\ms

\noindent{\bf Theorem \ref{Vogel}.} {\em Let $v_i$ denote the $i$th  St\"uckrad-Vogel cycle of the intersection of two equidimensional subschemes $X$ and $Y$ of ${\Bbb P}_k^n$, $i = 1,\ldots,n+1$. Let $I_X$ and $I_Y$ denote the defining ideals of $X$ and $Y$ in $k[x_0,\ldots,x_n]$ and $k[y_0,\ldots,y_n]$, respectively. Put $A = k[x_0,\ldots,x_n,y_0,\ldots,y_n]/(I_X,I_Y)$. Let $\frak m$ be the maximal graded ideal of $A$ and $J = (x_0-y_0,\ldots,x_n-y_n)A$. Then}
$$\deg v_i = e_{i-1}({\frak m}|J)-e_i({\frak m}|J).$$

This proposition together with a recent result of Achilles and Manaresi [AM] which interprets the degree of the St\"uckrad-Vogel cycles as mixed multiplicities of certain bigraded algebra (see Section 4) provide an interesting relationship between mixed multiplicities and intersection theory.\sk

This paper is divided into four sections. Sections 1 and Section 2 investigate the degrees and the mixed multiplicities of the Hilbert polynomial of a bigraded algebra, while Section 3 and Section 4 deal with the positivity of mixed multiplicities of ideals and with their applications. The notations introduced above will be kept throughout this paper. For unexplained terminologies we refer to Eisenbud's book on Commutative Algebra [E]. \sk

\noindent{\it Acknowledgment.} Corollary \ref{rigid3} has been also obtained by D. Katz (private communication).\sk 

\section{Partial degrees of the Hilbert polymomial}\sk

We shall need Northcott and Rees' theory on reductions of ideals [NR].\sk

Let $(A,\frak m)$ be a local ring and $I$ an arbitrary ideal of $A$. We call an ideal $J \subset I$ a {\it reduction} of $I$ if there is an integer $n \ge 0$ such that $I^{n+1} = JI^n$. A reduction of $I$ is said to be {\it minimal} if it does not contain any other reduction of $I$. These notions are closely related to the {\it fiber ring} of $I$ which is defined as the graded algebra
$$F(I) := \oplus_{n\ge 0}I^n/{\frak m}I^n.$$
In fact, it can be shown that $J$ is a reduction of $I$ if and only if the ideal of $F(I)$ generated by the degree one initial elements of $J$ is a primary ideal of the maximal graded ideal of $F(I)$. From this it follows that if the residue field of $A$ is infinite, the minimal number of generators of any minimal reduction $J$ of $I$ is equal to $\dim F(I)$. For this reason, $\dim F(I)$ is termed the {\it analytic spread} of $I$ and denoted by $s(I)$. We refer the reader to [NR] for more details. \sk

It is easy to verify that the results of [NR] also hold for homogeneous ideals generated by elements of the same degree in a standard graded algebra over an artinian local ring. Here we are interested only in the case $A = R$, where $R$ is a standard bigraded algebra over an artinian local ring $K$ and the ${\Bbb N}$-graded structure is given by $R_n = \oplus_{u+v=n}R_{(u,v)}$. \sk

Let $M$ denote the maximal homogeneous ideal of $R$. Notice that $M = {\frak n}R+R_{(1)}+R_{(2)}$, where $\frak n$ is the maximal ideal of $K$. 

\begin{Proposition}\label{spread}  $s(R_{(1)}) = \dim R/R_{(2)}$. \end{Proposition}

\begin{pf}  Since $R$ is standard bigraded,  $R_{(1)} = \oplus_{(u,v) \ge (1,0)}R_{(u,v)}$ and $R_{(2)} = \oplus_{(u,v) \ge (0,1)}R_{(u,v)}$, where $(u,v) \ge (a,b)$ means $u \ge a$ and $v \ge b$. Thus, $R_{(1)}^n = \oplus_{(u,v) \ge (n,0)}R_{(u,v)}$ and $MR_{(1)}^n = \oplus_{(u,v) \ge  (n+1,0)}R_{(u,v)}+ \oplus_{(u,v) \ge (n,1)}R_{(u,v)}+ {\frak n}R_{(n,0)}$. From this it follows that
$$F(R_{(1)}) = \oplus_{n \ge 0}R_{(1)}^n/MR_{(1)}^n  = \oplus_{n \ge 0}R_{(n,0)}/{\frak n}R_{(n,0)}.$$
Since $R/R_{(2)} = \oplus_{n \ge 0}R_{(n,0)}$, we get $F(R_{(1)}) = R/(R_{(2)}+{\frak n}R)$. Note that $\frak n$ is a nilpotent ideal. Then \sk
\hspace{2.5cm} $s(R_{(1)}) = \dim F(R_{(1)}) = \dim R/(R_{(2)}+{\frak n}R) = \dim R/R_{(2)}.$ \end{pf}

It is clear that any homogeneous minimal reduction of $R_{(1)}$ is generated by homogeneous elements of degree $(1,0)$. Now we shall see that these elements can be chosen to be like a regular sequence in large degree. \sk

A sequence $z_1,\ldots,z_s$ of homogeneous elements in $R$ is called {\it filter-regular} (with respect to $R_{++}$) if $z_i \not\in \wp$ for all associated prime ideals $\wp \not\supseteq R_{++}$ of $(z_1,\ldots,z_{i-1})$, $i = 1,\ldots,s$ (see e.g. [SV2] and [Tr] for the origin and basic properties).  We do not require that $z_i \not\in (z_1,\ldots,z_{i-1})$. Therefore, if $R_{++}$ is contained in the radical of $(z_1,\ldots,z_{i-1})$, any homogeneous element $z_i$ of $R$ satisfies the above condition. Note that for the above definition and the following observations we do not need to assume that $R_0$ is artinian.\sk

\begin{Lemma}\label{char}  A sequence $z_1,\ldots,z_s$ of homogeneous elements of $R$ is filter-regular  if and only if for $u$ and $v$ large enough,
$$[(z_1,\ldots,z_{i-1}):z_i]_{(u,v)} = (z_1,\ldots,z_{i-1})_{(u,v)},\ i = 1,\ldots,s.$$\end{Lemma}

\begin{pf} It is clear that $z_i \not\in \wp$ for all associated prime ideals $\wp \not\supseteq R_{++}$ of $(z_1,\ldots,z_{i-1})$ if and only if $[(z_1,\ldots,z_{i-1}):z_i]_{\wp'} = (z_1,\ldots,z_{i-1})_{\wp'}$ for all prime ideals $\wp' \not\supseteq R_{++}$. If we put 
$ E_i = (z_1,\ldots,z_{i-1}):z_i/(z_1,\ldots,z_{i-1})$, then the latter condition is equivalent to the condition that $E_i$ has only associated prime ideals which contain $R_{++}$. That means $E_i$ is annihilated by some power $R_{++}^n$. For $u$ and $v$ large enough we have
$(E_i)_{(u,v)} \subseteq R_{++}^nE_i.$ Hence $E_i$ is annihilated by some power $R_{++}^n$ if and only if  $(E_i)_{(u,v)} = 0$ or, equivalently, $[(z_1,\ldots,z_{i-1}):z_i]_{(u,v)} = (z_1,\ldots,z_{i-1})_{(u,v)}$ 
for $u$ and $v$ large enough. \end{pf}

\begin{Lemma}\label{exist} Let $Q$ be any reduction of $R_{(1)}$ generated by homogeneous elements of degree $(1,0)$. If the residue field of $K$ is infinite, there is a filter-regular sequence of homogeneous elements of degree $(1,0)$  which minimally generate $Q$. \end{Lemma}

\begin{pf} Since $Q$ and $R_{(1)}$ share the same radical and since $Q$ is generated by the elements of $Q_{(1,0)}$, any prime ideal $\wp$ of $R$ which does not contain $R_{(1)}$ also does not contain $Q_{(1,0)}$. By Nakayama's lemma, the set $\wp \cap Q_{(1,0)}$ maps to a proper subspace of the vector space $ Q_{(1,0)}/{\frak n}Q_{(1,0)}$ over the residue field $K/\frak n$. Since $K/\frak n$ is infinite, we may choose an element $z_1 \in Q_{(1,0)}\setminus {\frak n}Q_{(1,0)}$ such that $z_1 \not\in \wp$ for all associated prime ideals $\wp \not\supseteq Q_{(1,0)}$ of $R$. Since $R_{++} \subset R_{(1)}$, $z_1 \not\in \wp$ for any associated prime ideal $\wp \not\supseteq R_{++}$ of $R$.
Let $s$ be the  minimal number of generators of $Q$. It is clear that $Q/(z_1)$ is minimally generated by $s-1$ homogeneous elements of degree $(1,0)$. If $s > 1$, we can find, similarly as above, a homogeneous element $z_2$ of degree $(1,0)$ such that $z_1,z_2$ is a filter-regular sequence and $Q/(z_1,z_2)$ is minimally generated by $s-2$ elements of homogeneous elements of degree $(1,0)$. Continuing in this way we will find a filter-regular sequence $z_1,\ldots,z_s$  of homogeneous elements of degree $(1,0)$ such that $Q = (z_1,\ldots,z_s)$. \end{pf}

\noindent{\bf Remark.}  The assumption on the infiniteness of the residue field of $K$ does not cause us any problem because the Hilbert function does not change if we replace $R$ by $R \otimes_KK[t]_{{\frak n}K[t]}$, where $t$ is an indeterminate.\ms 

There is the following relationship between a filter-regular element of degree $(1,0)$ and  the Hilbert polynomial $P_R(u,v)$.

\begin{Proposition}\label{reduc1} Let $z$ be a homogeneous element of degree $(1,0)$ which is filter-regular. Assume that $P_R(u,v) \neq 0$. Then \par
{\rm (i) } $\deg_uP_{R/zR}(u,v) = \deg_uP_R(u,v)-1$, \par 
{\rm (ii) } $\deg P_{R/zR}(u,v) \le \deg P_R(u,v) -1$, \par
{\rm (iii) } $e_{ij}(R/zR) = e_{i+1j}(R)$ if $i+j \ge \deg P_R(u,v)-1$.
\end{Proposition}

\begin{pf} By Lemma \ref{char}, $(0:z)_{(u,v)} = 0$ for $u$ and $v$ large enough. This implies $P_{R/0:z}(u,v) = P_R(u,v).$ Consider the exact sequence $$0 \To R/0:z \overset z \To R \To R/zR \to 0.$$ Since $z$ is a homogeneous element of degree $(1,0)$, we have
$$\displaylines{P_{R/zR}(u,v) = P_R(u,v) - P_{R/0:z}(u-1,v) = P_R(u,v) - P_R(u-1,v) \cr
= \sum_{i+j \le r}e_{ij}(R) {u \choose i}{v \choose j}\ - \sum_{i+j \le r}e_{ij}(R){u-1 \choose i}{v \choose j}\cr
= \sum_{\tiny \begin{array}{c}i+j \le r \\ i \ge 1\end{array}} \left[{e_{ij}(R) \over (i-1)!j!}u^{i-1}v^j + \text{terms of degree $< i+j-1$ with $\deg_u < i-1$}\right],}$$
where $r = \deg P_R(u,v) \ge 0$ because $P_R(u,v) \neq 0$. Now we can easily deduce the statements of Proposition \ref{reduc1} from the above expression for $P_{R/zR}(u,v)$. \end{pf}

\noindent{\bf Remark.} If $P_R(u,v) = 0$, then $P_{R/zR}(u,v) = 0$. In this case, the degrees of $P_R(u,v)$ and $P_{R/zR}(u,v)$ can be any number.\ms

Proposition \ref{reduc1} leads to the following bound for the degree $\deg_uP_R(u,v)$ of $P_R(u,v)$ with respect to $u$.

\begin{Corollary}\label{less} $\deg_uP_R(u,v) \le \dim R/R_{(2)}-1$. \end{Corollary}

\begin{pf} Without restriction we may assume that the residue field of $K$ is infinite. Let $s = \dim R/R_{(2)}$. By Proposition \ref{spread} and Lemma \ref{exist} we can find a filter-regular sequence $z_1,\ldots,z_s$ of homogeneous elements of degree $(1,0)$ such that the ideal $Q = (z_1,\ldots,z_s)$ is a minimal reduction of $R_{(1)}$. Since there is an integer $n \ge 0$ such that $R_{(1)}^{n+1} = QR_{(1)}^n$, $(R/Q)_{(u,v)} = 0$ for $u > n$. Hence $P_{R/Q}(u,v) = 0$. So $\deg_uP_R(u,v) = -1$. By Proposition \ref{reduc1} (i) we must have\sk
\hspace{4cm} $\deg_uP_R(u,v) \le \deg_uP_{R/Q}(u,v) +s = s-1.$ \end{pf}
 
The above bound for $\deg_uP_R(u,v)$  helps us to establish the right formulae for the partial degrees of $P_R(u,v)$. To see this we shall need the following observations on the notation 
$${\frak a}:{\frak b}^\infty := \{x \in S|\ \text{there is a positive integer $n$ such that $x{\frak b}^n \subseteq {\frak a}$}\},$$
where $\frak a, \frak b$ are arbitrary ideals of a noetherian commutative ring $S$. First, there always exists an integer $n$ such that ${\frak b}^n({\frak a}:{\frak b}^\infty) \subseteq \frak a$ and, second, ${\frak a}:{\frak b}^\infty$ is the intersection of the primary components of $\frak a$ whose associated prime ideals do not contain $\frak b$. 

\begin{Lemma}\label{reduc2} $P_R(u,v) = P_{R/0:R_{++}^\infty}(u,v)$. \end{Lemma}

\begin{pf} Since $0:R_{++}^\infty$ is annihilated by some power of $R_{++}$, $(0:R_{++}^\infty)_{(u,v)} = 0$ for $u$ and $v$ large enough. Therefore, $R_{(u,v)} = (R/0:R_{++}^\infty)_{(u,v)}$ for $u$ and $v$ large enough. Hence $P_R(u,v) = P_{R/0:R_{++}^\infty}(u,v).$ \end{pf}

\begin{Theorem}\label{degree} Let $R$ be a standard bigraded algebra over an artinian local ring. Let $R_{(1)}$, $R_{(2)}$, and $R_{++}$ denote the ideals generated by the homogeneous elements of degree $(1,0)$, $(0,1)$, and $(1,1)$, respectively. Then
\begin{eqnarray*} \deg P_R(u,v) & = & \dim R/0:R_{++}^\infty-2,\\ \deg_uP_R(u,v) & = & \dim R/(0:R_{++}^\infty+R_{(2)})-1,\\ \deg_vP_R(u,v) & = & \dim R/(0:R_{++}^\infty+R_{(1)})-1. \end{eqnarray*}\end{Theorem}

\begin{pf} By Lemma \ref{reduc2} we have 
$$\deg P_R(u,v) = \deg P_{R/0:R_{++}^\infty}(u,v) = \rdim R/0:R_{++}^\infty.$$
Since the associated prime ideals of $0:R_{++}^\infty$ are exactly those of the zeroideal of $R$ which do not contain $R_{++}$, 
$$\rdim R/0:R_{++}^\infty = \dim R/0:R_{++}^\infty.$$
So we have proved the first formula of Theorem \ref{degree}.\par
To prove the second formula we have, by Lemma \ref{reduc2} and Corollary \ref{less},
$$\deg_u P_R(u,v) = \deg_uP_{R/0:R_{++}^\infty}(u,v) \le \dim R/(0:R_{++}^\infty+R_{(2)})-1.$$
Now we will show that $\deg_uP_R(u,v) \ge \dim R/(0:R_{++}^\infty+R_{(2)})-1.$ Let $a$ and $b$ be positive integers such that $h_R(u,v) = P_R(u,v)$ for $(u,v) \ge (a,b)$. Then
$$\deg_uP_R(u,v) \ge \deg_u P_R(u,b).$$
Since $R/R_{(2)} = \oplus_{u \ge 0}R_{(u,0)}$ is a standard ${\Bbb N}$-graded $K$-algebra and since $R_{(2)}^b/R_{(2)}^{b+1} = \oplus_{u \ge 0}R_{(u,b)}$, $P_R(u,b)$ is the Hilbert polynomial of the graded $(R/R_{(2)})$-module $R_{(2)}^b/R_{(2)}^{b+1}$. Therefore,
$$\deg_u P_R(u,b) = \dim R_{(2)}^b/R_{(2)}^{b+1}-1.$$
It remains to show that $\dim R_{(2)}^b/R_{(2)}^{b+1} \ge \dim R/(0:R_{++}^\infty+R_{(2)})$. Let $\wp$ be an associated prime ideal of $0:R_{++}^\infty+R_{(2)}$ such that $\dim R/\wp = \dim R/(0:R_{++}^\infty+R_{(2)})$. Note that $0:R_{(2)}^b \subseteq 0:R_{++}^\infty \subseteq \wp$. Then $(0:R_{(2)}^b)_\wp$ is a proper ideal of $R_\wp$, whence $(R_{(2)}^b)_\wp \neq 0$. By Nakayama's lemma, this implies $(R_{(2)}^b)_\wp \neq (R_{(2)}^{b+1})_\wp$. Hence $(R_{(2)}^b/R_{(2)}^{b+1})_\wp \neq 0$. Thus,
$$\dim R_{(2)}^b/R_{(2)}^{b+1} \ge \dim R/\wp = \dim R/(0:R_{++}^\infty+R_{(2)}).$$
So we have proved the second formula of Theorem \ref{degree}. The third formula can be proved similarly by changing the order of the bidegree of $R$. \end{pf}

It is quite natural to ask whether $\deg_uP_R(u,v) = \dim R/R_{(2)}-1$ and $\deg_vP_R(u,v) = \dim R/R_{(1)}-1$. The following example shows that we may have $P_R(u,v) = 0$ while $R, R/R_{(1)}$ and $R/R_{(2)}$ have arbitrary dimension.\ms

\noindent{\bf Example.} Let $R = K[x_1,\ldots,x_n,y_1,\ldots,y_n]/(x_1,\dots,x_n)\cap(y_1,\dots,y_n)$ with $\deg x_i =  (1,0)$, $\deg y_i = (0,1)$, $i = 1,\ldots,n$. Then $R_{++} = 0$. Therefore $P_R(u,v) = 0$ while $\dim R = \dim R/R_{(1)} = \dim R/R_{(2)} = n$. \ms

To compute $\deg_uP_R(u,v)$ and $\deg_vP_R(u,v)$ we will also use the following equalities.

\begin{Lemma} \label{degree2} Let $R_{++}$, $R_{(1)}$ and $R_{(2)}$ be as above. Then
\begin{eqnarray*} \dim R/(0:R_{++}^\infty+R_{(1)}) & = & \dim R/(0:R_{(1)}^\infty+R_{(1)}),\\ \dim R/(0:R_{++}^\infty+R_{(2)}) & = & \dim R/(0:R_{(2)}^\infty+R_{(2)}). \end{eqnarray*}
\end{Lemma}

\begin{pf} Since $R_{++} = R_{(1)}R_{(2)}$, there is an integer $n$ such that 
$$(R_{(1)}+R_{(2)})^n(0:R_{++}^\infty + R_{(1)}) \subseteq  0:R_{(1)}^\infty + R_{(1)}  \subseteq 0:R_{++}^\infty + R_{(1)}.$$ 
Since $\sqrt{R_{(1)}+R_{(2)}}$ is the maximal graded ideal of $R$, this implies that 
$$\sqrt{0:R_{++}^\infty + R_{(1)}} = \sqrt{0:R_{(1)}^\infty + R_{(1)}}.$$
Therefore, $\dim R/(0:R_{++}^\infty+R_{(1)}) = \dim R/(0:R_{(1)}^\infty+R_{(1)})$. 
Similarly, we also get $\dim R/(0:R_{++}^\infty+R_{(2)}) = \dim R/(0:R_{(2)}^\infty+R_{(2)})$.\end{pf}

\section{Mixed multiplicities of bigraded algebras}\sk

Let $R$ be an arbitrary standard bigraded algebra over an artinian local ring $K$. For simplicity we set 
\begin{eqnarray*} r & := & \dim R/0:R_{++}^\infty-2,\\ r_1 & := & \dim R/(0:R_{++}^\infty+R_{(2)})-1,\\ r_2 & := & \dim R/(0:R_{++}^\infty+R_{(1)})-1. \end{eqnarray*}
By Theorem \ref{degree} we have $ \deg P_R(u,v) = r,\ \deg_uP_R(u,v) = r_1,\ \deg_vP_R(u,v) = r_2.$ Therefore $r_1 + r_2 \ge r$ and we have the following estimate for the numbers $e_{ij}(R)$ of the Hilbert polynomial $P_R(u,v)$. 

\begin{Lemma}\label{zero} $e_{ij}(R) = 0$ for $i > r_1$ or $j > r_2$ or $i+j > r$. \end{Lemma}

Now we will want to study the mixed multiplicities $e_{ij}(R)$, $i+j = r$. For this we shall need the following relationship between Samuel's multiplicity and mixed multiplicities.\sk

Let $R_n = \oplus_{u+v = n}R_{(u,v)}$. Then $R = \oplus_{n\ge 0}R_n$ is a standard graded algebra over $K$. It is well-known that there is a polynomial $P_R(n)$ of degree $d-1$, $d = \dim R$, such that $\ell(R_n) = P_R(n)$ for $n$ large enough. Write 
$$P_R(n) = {a \over (d-1)!}n^{d-1} + \text{terms of degree $< d-1$}.$$
Then $e(R) = a$ is the  multiplicity of the graded algebra $R$. 

\begin{Proposition} \label{mult} {\rm [KMV, Theorem 4.1]} Assume that $\height R_{(1)} \ge 1$ and $\height R_{(2)} \ge 1$. Then 
$$e(R) = \sum_{i+j=d-2} e_{ij}(R).$$ \end{Proposition}

Moreover, we shall need the existence of certain filter-regular sequences.

\begin{Lemma}\label{exist2} Let $i$ and $j$ be non-negative integers, $i+j = r$. If the residue field of $K$ is infinite, there exist homogeneous elements $x_1,\ldots,x_i$ and $y_1,\ldots,y_j$ of degrees $(1,0)$ and $(0,1)$, respectively, such that $x_1,\ldots,x_i,y_1,\ldots,y_j$ is a filter-regular sequence. \end{Lemma}

\begin{pf} By Lemma \ref{exist} we can find a filter-sequences $x_1,\ldots,x_i$ of homogeneous elements of degree $(1,0)$. Note that if $i$ is greater than the minimal number of generators of $R_{(1)}$, we may add arbitrary homogeneous elements to an existing filter-regular sequence which generates $R_{(1)}$ to get a new filter-regular sequence. Similarly, we can also find homogeneous elements $y_1,\ldots,y_j$ of degree $(0,1)$ in $R$ which form a filter-regular sequence in $R/(x_1,\ldots,x_i)$. It is clear from the definition of filter-regular sequences that $x_1,\ldots,x_i,y_1,\ldots,y_j$ is a filter-regular sequence in $R$. \end{pf}

The following result gives an effective criterion for the positivity of a mixed multiplicity $e_{ij}(R)$ and shows how to compute $e_{ij}(R)$ as the multiplicity of a graded algebra.

\begin{Theorem} \label{nonzero} Let $i, j$ be non-negative integers, $i+j = r$. Let $x_1,\ldots,x_i$ be a filter-regular sequence of homogeneous elements of degree $(1,0)$. Then $e_{ij}(R) > 0$ if and only if 
$$\dim R/((x_1,\ldots,x_i):R_{++}^\infty+R_{(1)}) = j+1.$$ In this case, if we choose homogeneous elements $y_1,\ldots,y_j$ of degree $(0,1)$ such that $x_1,\ldots,x_i,y_1,\ldots,y_j$ is a filter-regular sequence, then  
$$e_{ij}(R) = e(R/(x_1,\ldots,x_i,y_1,\ldots,y_j):R_{++}^\infty).$$ \end{Theorem}

\begin{pf} Let $Q_h = (x_1,\ldots,x_h)$,   $h = 1,\ldots,i$. Applying Proposition \ref{reduc1} successively to the bigraded algebras $R/Q_h$ we will obtain 
$$\displaylines{\deg P_{R/Q_i}(u,v) \le d-2-i = j,\cr e_{ij}(R) = e_{0j}(R/Q_i).}$$ 
If $\dim R/(Q_i:R_{++}^\infty+R_{(1)}) > j+1$, then $\deg P_{R/Q_i}(u,v) \ge \deg_vP_{R/Q_i}(u,v) > j$ by Theorem \ref{degree}.  Hence we obtain a contradiction. If $\dim R/(Q_i:R_{++}^\infty+R_{(1)}) = j+1$, then $\deg P_{R/Q_i}(u,v) \ge \deg_vP_{R/Q_i}(u,v) = j$ by Theorem \ref{degree}.
Therefore, we must have $\deg P_{R/Q_i}(u,v) = \deg_vP_{R/Q_i}(u,v) = j.$ This means that the coefficient of the monomial $v^j$ in $P_{R/Q_i}(u,v)$ is positive. Hence  $e_{ij}(R) = e_{0j}(R/Q_i) > 0$. If $\dim R/(Q_i:R_{++}^\infty+R_{(1)}) < j+1$, then $e_{ij}(R) = e_{0j}(R/Q_i) = 0$ by Lemma \ref{zero}. So we get the first statement of Theorem \ref{nonzero}. \par

To prove the second statement we apply Proposition \ref{reduc1} successively to the bigraded algebras $ R/(Q_i,y_1,\ldots,y_h)$, $h = 1,\ldots,j$. If we put $Q=(x_1,\ldots,x_i,y_1,\ldots,y_j)$, then 
$$\displaylines{\deg P_{R/Q}(u,v) \le \deg P_{R/Q_i}(u,v) - j = 0\cr
e_{00}(R/Q) = e_{0j}(R/Q_i) = e_{ij}(R) > 0.}$$
>From this it follows that $\deg P_{R/Q}(u,v) = 0.$ Let 
$$\bar R = R/Q:R_{++}^\infty =  R/(x_1,\ldots,x_i,y_1,\ldots,y_j):R_{++}^\infty.$$  
Note that $0:\bar R_{++}^\infty = 0$.  Since $0:\bar R_{(1)}$ and $0:\bar R_{(2)}$ are contained in $0:\bar R_{++}^\infty$,  $0:\bar R_{(1)} = 0:\bar R_{(2)} = 0$. Therefore,  $\height \bar R_{(1)} \ge 1$ and $\height \bar R_{(2)} \ge 1$. By Lemma \ref{reduc2}, $P_{\bar R}(u,v) = P_{R/Q}(u,v)$. Hence, using Theorem \ref{degree} we get  $\dim \bar R = \deg P_{\bar R}(u,v)+2 = \deg P_{R/Q}(u,v)+2 = 2$. Now we apply  Proposition \ref{mult} and obtain 
$$e(\bar R) =  e_{00}(\bar R) = e_{00}(R/Q) = e_{ij}(R).$$
This completes the proof of Theorem \ref{nonzero}. \end{pf}

By Lemma \ref{zero} we know that $e_{ir-i}(R) = 0$ for $i < r-r_2$ or $i > r_1$. Using Theorem \ref{nonzero} we can construct examples with $e_{ir-i}(R) = 0$ for $i = r-r_2,r_1$. \ms

\noindent{\bf Example.} Let $R = K[X,Y]/(x_1,y_1) \cap (x_1,x_2,x_3) \cap (y_1,y_2,y_3)$ with $X = \{x_1,x_2,x_3,x_4\}$, $Y = \{y_1,y_2,y_3,y_4\}$ and $\deg x_i = (1,0),\ \deg y_i = (0,1),\ i = 1,2,3,4$. Then $R/R_{(1)} \cong K[Y]$ and $R/R_{(2)} \cong K[X]$. Since $0:R_{++}^\infty = 0$, we get 
$$r = \dim R - 2 = 4,\ r_1 = \dim R/R_{(1)} - 1= 3,\ r_2 = \dim R/R_{(2)} - 1 = 3.$$
It is clear that $x_4$ is a non-zerodivisor in $R$ and  
$$x_4R:R_{++}^\infty = (x_1,y_1,x_4)R \cap (y_1,y_2,y_3,x_4)R.$$ 
Then $R/(x_4R:R_{++}^\infty + R_{(1)}) \cong K[Y]/(y_1)$. Since $\dim R/(x_4R:R_{++}^\infty + R_{(1)}) = 3 < 3+1$, $e_{13}(R) = 0$. By symmetry we also have $e_{31}(R) = 0$. Now we want to compute the only non-vanishing mixed multiplicity $e_{22}(R)$ of $R$. It is easy to check that $x_4,x_2,y_4,y_2$ is a filter-regular sequence in $R$.
Put $Q = (x_4,x_2,y_4,y_2)$. Then
$$R/Q = K[X,Y]/(x_1,x_2,x_4,y_1,y_2,y_4) \cap (x_1,x_2,x_3,x_4,y_2,y_4) \cap (x_2,x_4,y_1,y_2,y_3,y_4).$$
>From this it follows that $R/Q:R_{++}^\infty = K[X,Y]/(x_1,x_2,x_4,y_1,y_2,y_4) \cong K[x_3,y_3].$ Hence $e_{22}(R) = e(R/Q:R_{++}^\infty) = \ell(K).$\ms

\begin{Corollary} \label{first}  $e_{0r}(R) > 0$ $(e_{r0}(R) > 0)$ if and only if $r_2+1 = r$ $(r_1+1 = r)$. \end{Corollary} 

\begin{pf}  Since $\dim R/(0:R_{++}^\infty+R_{(1)}) = r_2+1$ ($\dim R/(0:R_{++}^\infty+R_{(2)}) = r_1+1$), the statement follows from Theorem \ref{nonzero} (by changing the order of the bidegree of $R$, respectively). \end{pf}

The sequence $e_{r0}(R),e_{r-11}(R),\ldots,e_{0r}(R)$ of mixed multiplicities may behave wildly. Katz, Mandal and Verma [KMV, Example 5.2] showed that it can be any sequence of non-negative integers $c_0,c_1,\ldots,c_r$ with at least a positive entry. They also raised the question whether it is {\it rigid}, i.e. there are integers $a, b$ such that $e_{ir-i}(R) > 0$ for $a \le i \le b$ and $e_{ir-i}(R) = 0$ otherwise,  if $R$ is a domain or Cohen-Macaulay [KMV, Question 5.3]. We shall see that this question has a positive answer by showing that $e_{ir-i}(R) > 0$ for $r-r_2 \le i \le r_1$.\sk

We  need the following version of Grothendieck's Connectedness Theorem which is due to Brodmann and Rung [BR, Proposition 2.1]. For any commutative noetherian ring $S$ let
\begin{eqnarray*} c(S) & := & \min\{\dim Z|\ Z \subseteq \Spec(S)\ \text{is closed and $\Spec(S)\setminus Z$ is disconnected}\},\\
\sdim S & := & \min\{\dim S/\wp|\ \wp\ \text{is a minimal associated prime ideal of $S$}\}. \end{eqnarray*}

\begin{Lemma} \label{connect} Let $P$ and $Q$ be proper ideals of $R$ such that 
$$\min\{\dim R/P, \dim R/Q\} > \dim R/(P+Q).$$ Let $\ara(P \cap Q)$ denote the minimum number of generators of ideals which have the same radical as $P \cap Q$. Then
$$\dim R/(P+Q) \ge \min\{c(R), \sdim R -1\} - \ara(P\cap Q).$$ \end{Lemma}

\begin{pf} The statement was originally proved for ideals of a complete local ring [BS, Lemma 19.2.8]. As noted in [BR, Section 1], the proof can be easily extended to homogeneous ideals of a standard graded ring over a local complete ring. Note that an artinian local ring is always complete. \end{pf}
	 
Recall that a commutative noetherian ring $S$ is said to satisfy the {\it first chain condition} if maximal chains of prime ideals in $S$ have the same length. Moreover, we say that $S$ is {\it connected in codimension} 1 if $c(S) = \dim S-1$.

\begin{Theorem} \label{rigid1} Assume that $R/0:R_{++}^\infty$ satisfies the first chain condition. Then $e_{ir-i}(R)> 0$ for $i = r-r_2, r_1$. If $R/0:R_{++}^\infty$ is moreover connected in codimension 1, then  $e_{ir-i}(R)> 0$ for $r-r_2 \le i \le r_1$. \end{Theorem}

\begin{pf} By Lemma \ref{reduc2} we may replace $R$ by $R/0:R_{++}^\infty$. Then $0:R_{++}^\infty = 0$. Hence $\dim R = r+2$. Without restriction we may assume that the residue field of $K$ is infinite. By Lemma \ref{spread}, $r_1+1 =\dim R/R_{(2)} =  s(R_{(1)})$ is the minimal number of generators of minimal reductions of $R_{(1)}$. Therefore, using Lemma \ref{exist} we can find a filter-regular sequence $x_1,\ldots,x_{r_1+1}$ of homogeneous elements of degree $(1,0)$ such that $(x_1,\ldots,x_{r_1+1})$ is a reduction of $R_{(1)}$.
Let $Q_i = (x_1,\ldots,x_i)$,  $r-r_2 \le i \le r_1$. We will first estimate $\dim R/Q_i:R_{++}^\infty$. Let $\wp$ be a minimal associated prime ideal of $Q_i:R_{++}^\infty$ with $\dim R/\wp = \dim R/Q_i:R_{++}^\infty$. Then $\wp \not\supseteq R_{++}$. By the definition of filter-regular sequences, $x_h$ does not belong to any associated prime ideal of $Q_{h-1}R_\wp$, $h = 1,\ldots,i$. Therefore, $x_1,\ldots,x_i$ is a regular sequence in $R_\wp$. From this it follows that $\dim R_\wp = i$. Since $R$ satisfies the first chain condition, $\dim R/\wp = \dim R - \height \wp = r+2 - \dim R_\wp = r-i+2$ so that we get
$$\dim R/Q_i:R_{++}^\infty = r - i +2.$$
Since $x_{i+1}$  does not belong to any associated prime ideals of $Q_i:R_{++}^\infty$,  
$$\dim R/(Q_i:R_{++}^\infty + x_{i+1}R) = \dim R/Q_i:R_{++}^\infty -1 = r-i+1.$$
In particular, $\dim R/(Q_{r_1}:R_{++}^\infty+x_{r_1+1}R) = r-r_1+1$. Since $R_{(1)}$ and $(x_1,\ldots,x_{r_1+1}) = Q_{r_1}+x_{r_1+1}R$ share the same radical, we also have $\dim R/(Q_{r_1}:R_{++}^\infty + R_{(1)}) = r-i+1$. By Theorem \ref{nonzero}, this implies $e_{r_1r-r_1}(R)> 0$. By changing the order of the bidegree of $R$ we can also show that $e_{r-r_2r_2}(R) > 0$.\par  
Now let $r-r_2 < i < r_1$. Since $x_{i+1} \in R_{(1)}$,
$$\dim R/(Q_i:R_{++}^\infty + R_{(1)}) \le \dim R/(Q_i:R_{++}^\infty + x_{i+1}R) = r-i+1.$$
Applying Lemma \ref{degree2} to the bigraded algebra $R/Q_i$ we have 
$$ \dim R/(Q_i:R_{(1)}^\infty + R_{(1)}) = \dim R/(Q_i:R_{++}^\infty + R_{(1)}) \le r-i+1.$$
Write $Q_i = (Q_i:R_{(1)}^\infty) \cap Q$ where $Q$ is an ideal whose associated prime ideals contain $R_{(1)}$. Since $i < r_1 = s(R_{(1)})-1$, $J_i$ is not a reduction of $R_{(1)}$. Since $R$ is standard bigraded, this implies that $R_{(1)}$ is not contained in the radical of $J_i$. Hence $Q_i:R_{(1)}^\infty$ is a proper ideal of $R$. Since $Q_i \subseteq R_{(1)}$, we may assume that $Q \subseteq R_{(1)}$. Then $\sqrt{Q} = \sqrt{R_{(1)}}$. Hence
$$\dim R/Q = \dim R/R_{(1)} = r_2 + 1 \ge r-i+2 > \dim R/(Q_i:R_{(1)}^\infty + R_{(1)})$$
Since $Q_i:R_{(1)}^\infty \subseteq Q_i:R_{++}^\infty$, we also have 
$$\dim R/Q_i:R_{(1)}^\infty \ge \dim R/Q_i:R_{++}^\infty = r-i+2 > \dim R/(Q_i:R_{(1)}^\infty + R_{(1)}).$$
Now we can apply Lemma \ref{connect} to the ideals $Q_i:R_{(1)}^\infty$ and $Q$ with $\ara((Q_i:R_{(1)}^\infty) \cap Q) = \ara(Q_i) \le  i$. If $R/0:R_{++}^\infty$ is connected in codimension 1,  $c(R) = \dim R -1 = r+1$. Since $\sdim R = \dim R = r+2$,  
$$\dim R/(Q_i:R_{(1)}^\infty+Q) \ge r+1 - \ara((Q_i:R_{(1)}^\infty) \cap Q) \ge r-i+1.$$
Since $R_{(1)} \subseteq \sqrt{Q}$,
\begin{eqnarray*} \dim R/(Q_i:R_{(1)}^\infty+R_{(1)}) & \ge & \dim R/\sqrt{Q_i:R_{(1)}^\infty + Q}\\
& = & \dim R/(Q_i:R_{(1)}^\infty+Q)\ \ge\ r-i+1. \end{eqnarray*} 
So we can conclude that
$$\dim R/(Q_i:R_{++}^\infty+R_{(1)}) = \dim R/(Q_i:R_{1)}^\infty+R_{(1)}) =  r-i+1.$$
By Theorem 2.4, this implies  $e_{ir-i}(R) > 0$. The proof of Theorem \ref{rigid1} is now complete. \end{pf}

\begin{Corollary} \label{rigid2} Let $R$ be a domain or a Cohen-Macaulay ring with $\height R_{++} \ge 1$. Then $e_{ir-i}(R)> 0$ for $r-r_2 \le i \le r_1$. \end{Corollary}

\begin{pf} If $R$ is a domain, $0:R_{++}^\infty = 0$ and $K$ must be a field. Hence $R$ is a factor ring of a polynomial ring over a field by a prime ideal. Since such a polynomial ring is catenary, that is, the maximal chains of primes ideals between two prime ideals $\wp \supset \wp'$ have the same length, $R$ satisfies the first chain condition. It is clear that  $c(R) = \dim R$ in this case. If $R$ is a Cohen-Macaulay ring, $0:R_{++}^\infty = 0$. It is well-known that $R$ satisfies the first chain condition. Moreover, $R$ is connected in codimension 1 by Hartshorne's Connectedness Theorem (see e.g. [E, Theorem 8.12]). Now we only need to apply Theorem \ref{rigid1} to both cases. \end{pf}

\section{Mixed multiplicities of ideals}\sk

Throughout this section we set $R = R(I|J)$, where $I$ is an $\frak m$-primary ideal, $J$ is an arbitrary  ideal of a local ring $(A,\frak m)$ and
$$R(I|J) := \oplus_{v \ge 0} I^vJ^u/I^{v+1}J^u.$$ 
By the definition of the mixed multiplicities of $I$ and $J$ we have $e_{i}(I|J) := e_{ir-i}(R),$ where $r = \deg P_R(u,v)$. \sk

If $\height J \ge 1$, Katz and Verma  [KV, Lemma 2.2] showed that $\deg P_R(u,v) = \dim A -1$ (implicitly) and $e_0(I|J) = e(I,A)$, where $e(I,A)$ denotes the Samuel's multiplicity of $A$ with respect to $I$. To estimate the other mixed multiplicities we need to reformulate their result for an ideal $J$ of arbitrary height.  

\begin{Lemma}\label{e0} For an arbitrary ideal $J$ of $A$ we have \par
{\rm (i) } $\deg P_R(u,v) = \dim A/0:J^\infty-1$,\par
{\rm (ii) } $e_0(I|J) = e(I,A/0:J^\infty)$. \end{Lemma}

\begin{pf} If $J$ is nilpotent, $0:J^\infty = A$, hence the statements are trivial. If $J$ is not a nilpotent ideal, we set $\bar A = A/0:J^\infty$, $\bar I = I\bar A$ and $\bar J = J\bar A$. By [KV, Lemma 2.3], 
$$P_R(u,v) = P_{R(\bar I|\bar J)}(u,v).$$
>From this it follows that
$$e_0(I|J) = e_0(\bar I|\bar J).$$
On the other hand, since $0:\bar J^\infty = 0$, $\bar J$ is not contained in any associated prime ideal of $\bar A$. This implies $\height \bar J \ge 1$. Applying [KV, Lemma 2.2] to the factor ring $\bar A$ we get 
$$ \deg P_{R(\bar I|\bar J)}(u,v) = \dim A/0:J^\infty-1$$
$\hspace{5.5cm}e_0(\bar I|\bar J) = e(I,A/0:J^\infty).$ \end{pf}

Katz and Verma also showed that $e_i(I|J) = 0$ for $i \ge s(J)$ [KV, Theorem 2.7 (i)]. But this is just a consequence of the following bound for the partial degree $\deg_uP_R(u,v)$ of $ P_R(u,v)$ with respect to the variable $u$.  

\begin{Proposition}\label{deg u} $\deg_u P_R(u,v) < s(J)$. \end{Proposition}

\begin{pf} By Lemma \ref{less}, $\deg_u P_R(u,v) < \dim R/R_{(2)}$. On the other hand, 
$$R/R_{(2)} = \oplus_{u \ge 0}R_{(u,0)} = \oplus_{u \ge 0}J^u/IJ^u = R(J)/IR(J).$$ 
Since $I$ is a $\frak m$-primary ideal, $\dim R(J)/IR(J) = \dim R(J)/{\frak m}R(J)$. But $R(J)/{\frak m}R(J) = \oplus_{u \ge 0}J^u/{\frak m}J^u = F(J).$ Therefore, $\dim R/R_{(2)} = \dim F(J) = s(J)$. \end{pf}
 
To reduce the computation of the mixed multiplicities $e_i(I|J)$ to the case $i = 0$ we need to consider the standard bigraded algebra
$$R(J|I) :=  \oplus_{(u,v) \in {{\Bbb N}}^2} I^vJ^u/I^vJ^{u+1}.$$

\begin{Lemma} \label{commute} Let $a_1,\ldots,a_i$, $0\le i < s(J)$, be elements in $J$ such that their images $x_1,\ldots,x_i$ in $J/IJ$ and $x'_1,\ldots,x'_i$ in $J/J^2$ form filter-regular sequences in $R(I|J)$ and $R(J|I)$, respectively. Let $\bar R = R(I|J)/(x_1,\ldots,x_i)$ and $\bar R' = R(J|I)/(x'_1,\ldots,x'_i)$. Let $\bar A = A/(a_1,\ldots,a_i)$, $\bar I = I\bar A$ and $\bar J = J\bar A$. Then there are surjective graded homomorphisms from $\bar R$ to $R(\bar I|\bar J)$ and from $\bar R'$ to $R(\bar J|\bar I)$ such that
\begin{eqnarray*} \bar R_{(u,v)} & = & R(\bar I|\bar J)_{(u,v)},\\
 \bar R'_{(u,v)} & = & R(\bar J|\bar I)_{(u,v)} \end{eqnarray*}
for $u$ and $v$ large enough. \end{Lemma}

\begin{pf} We first consider the case $i = 1$. For all $(u,v) \in {{\Bbb N}}^2$ we have
\begin{eqnarray*} \bar R_{(u,v)}& = & I^vJ^u/(I^{v+1}J^u + a_1I^vJ^{u-1}),\\
R(\bar I|\bar J)_{(u,v)} & = & \bar I^v\bar J^u/\bar I^{v+1}\bar J^u\\ &
= & (I^vJ^u + (a_1))/(I^{v+1}J^u + (a_1))\\
& \cong & I^vJ^u/(I^{v+1}J^u + (a_1)) \cap I^vJ^u\\
& = & I^vJ^u/(I^{v+1}J^u + (a_1)\cap I^vJ^u) \end{eqnarray*}
Therefore, there is a surjective graded homomorphism from $\bar R$ to 
$R(\bar I, \bar J)$. By the Artin-Rees lemma, there are integers $u_0$ and $v_0$ such that 
$$(a_1)\cap I^vJ^u = J^{u-u_0}((a_1) \cap I^uJ^{u_0}) = I^{u-u_0}J^{v-v_0}((a_1) \cap I^{v_0}J^{u_0}) \subseteq a_1 I^{u-u_0}J^{v-v_0}$$
and hence $I^vJ^u:a_1 \subseteq (0:a_1) + I^{v-v_0}J^{u-u_0}$ for $u \ge u_0$ and $v \ge v_0$. It follows that
$$(a_1) \cap I^vJ^u = a_1[(I^vJ^u:a_1) \cap I^{v-v_0}J^{u-u_0}].$$
By Lemma \ref{char}, $(0:x_1)_{(m,n)} = 0$ for $m$ and $n$ large enough. This can be translated as
$$(I^{n+1}J^{m+1}:a_1) \cap I^nJ^m = I^{n+1}J^m.$$
Using this relation for $n = u-u_0$ and $m = v-v_0,\ldots,v-1$ we get
\begin{eqnarray*} (I^vJ^u:a_1) \cap I^{v-v_0}J^{u-u_0}  & =  & (I^vJ^u:a_1) \cap (I^{v-v_0+1}J^{u-u_0+1}:a_1) \cap I^{v-v_0}J^{u-u_0}\\ &  = & (I^vJ^u:a_1) \cap I^{v-v_0+1}J^{u-u_0}  = \cdots \\ & = & (I^vJ^u:a_1) \cap (I^vJ^{u-u_0+1}:a_1) \cap I^{v-1}J^{u-u_0}\\ &  = & (I^vJ^u:a_1) \cap I^vJ^{u-u_0}. \end{eqnarray*}
On the other hand, by Lemma \ref{char}, $(0:x_1')_{(m,n)} = 0$ for $m$ and $n$ large enough. This can be translated as
$$(I^nJ^{m+2}:a_1) \cap I^nJ^m = I^nJ^{m+1}.$$ 
Using this relation for $m = u-u_0,\ldots,u-2$ and $n = v$ we can show similarly as above that for $u \ge u_0$ and $v \ge v_0$,  $$ (I^vJ^u:a_1) \cap I^vJ^{u-u_0} = \cdots = (I^vJ^u:a_1) \cap I^vJ^{u-2} = I^vJ^{u-1}.$$
Summing up we obtain $(a_1) \cap I^vJ^u = a_1I^vJ^{u-1}$ and therefore $\bar R _{(u,v)} = R(\bar I|\bar J)_{(u,v)}$ for $u$ and $v$ large enough. In the same way we can also show that there is a graded surjective homomorphism from $\bar R'$ to $R(\bar J|\bar I)$ such that  $\bar R' _{(u,v)} = R(\bar J|\bar I)_{(u,v)}$. \par
Now let $i > 1$ and $S = R(I/(a_1)|J/(a_1))$ and $S' = R(J/(a_1)|I/(a_1))$. As we have seen above, there are surjective graded homomorphisms from $R/(x_1)$ to $S$ and from $R'/(x'_1)$ to $S'$ such that $(R/(x_1)) _{(u,v)} = S_{(u,v)}$ and $(R'/(x'_1)) _{(u,v)} = S'_{(u,v)}$  for $u$ and $v$ large enough. They induce surjective graded homomorphisms from $\bar R$ to $S/(x_2,\ldots,x_i)$ and from $\bar R'$ to $S'/(x'_2,\ldots,x'_i)$ such that 
\begin{eqnarray*} \bar R_{(u,v)} & = & (S/(x_2,\ldots,x_i))_{(u,v)},\\
\bar R'_{(u,v)} & = & (S'/(x'_2,\ldots,x'_i))_{(u,v)},\end{eqnarray*}
for $u$ and $v$ large enough. Using Lemma \ref{char} we can show that $x_2,\ldots,x_i$ and $x'_2,\ldots,x'_i$ form filter-regular sequences in $S$ and $S'$, respectively. By induction we may asssume that there are surjective graded homomorphisms from $S/(x_2,\ldots,x_i)$ to $R(\bar I|\bar J)$ and from $S'/(x'_2,\ldots,x'_i)$ to $R(\bar J|\bar I)$ such that
\begin{eqnarray*} (S/(x_2,\ldots,x_i))_{(u,v)}&  = & R(\bar I|\bar J)_{(u,v)},\\
(S'/(x'_2,\ldots,x'_i))_{(u,v)}&  = & R(\bar J|\bar I)_{(u,v)},\end{eqnarray*}
for $u$ and $v$ large enough. So we can conclude that there are surjective graded  homomorphisms from $\bar R$ to $R(\bar I|\bar J)$ and from $\bar R'$ to $R(\bar J|\bar I)$ as stated in Lemma \ref{commute}. \end{pf}

Now we can say exactly when $e_i(I|J) > 0$ and how one can compute $e_i(I|J)$ for  $0 \le i < s(J)$.

\begin{Theorem} \label{main} Let $J$ be an arbitrary ideal of $A$ and $0 \le i < s(J)$. Let $a_1,\ldots,a_i$ be elements in $J$ such that their images in $J/IJ$ and $J/J^2$ form filter-regular sequences in $R(I|J)$ and $R(J|I)$, respectively. Then $e_i(I|J) > 0$ if and only if  $$\dim A/(a_1,\ldots,a_i):J^\infty = \dim A/0:J^\infty-i.$$ In this case, $e_i(I|J) = e(I,A/(a_1,\ldots,a_i):J^\infty).$  \end{Theorem}

\begin{pf}  Let $r = \dim A/0:J^\infty -2$. Then $\deg P_R(u,v) = r$ by Lemma \ref{e0}. Hence
$$e_i(I|J) = e_{ir-i}(R).$$ Let $x_1,\ldots,x_i$ be  the images of $a_1,\ldots,a_i$ in $J/IJ$. Put $\bar R = R/(x_1,\ldots,x_i)$. Applying Proposition \ref{reduc1} succesively to $R/(x_1,\ldots,x_h)$, $h = 1,\ldots,i$, we get 
$$\displaylines{\deg P_{\bar R}(u,v)  \le r-i,\cr
e_{ir-i}(R) = e_{0r-i}(\bar R).}$$
Let  $\bar A = A/(a_1,\ldots,a_i)$, $\bar I = I\bar A$ and $\bar J = J\bar A$. 
By Lemma \ref{commute}, $P_{\bar R}(u,v) = P_{R(\bar I|\bar J)}(u,v).$ Hence 
$$\displaylines{\deg P_{R(\bar I|\bar J)}(u,v) \le r-i,\cr
e_{0r-i}(\bar R)  = e_{0r-i}(R(\bar I|\bar J)).}$$
So we obtain $e_i(I|J) = e_{0r-i}(R(\bar I|\bar J)).$ According to Lemma \ref{e0} we have
$$\deg P_{R(\bar I|\bar J)}(u,v) = \dim \bar A/0:\bar J^\infty - 2 = \dim A/(x_1,\ldots,x_i):J^\infty - 2.$$ Hence $\dim A/(x_1,\ldots,x_i):J^\infty \le r-i+2$. \par
If $\dim A/(x_1,\ldots,x_i):J^\infty < r-i+2$, then $\deg P_{R(\bar I|\bar J)}(u,v) < r-i$. Therefore, $e_{0r-i}(R(\bar I|\bar J)) = 0$. If $\dim A/(x_1,\ldots,x_i):J^\infty = r-i+2$, then $\deg P_{R(\bar I|\bar J)}(u,v) = r-i$. Now we apply Lemma \ref{e0} again and obtain
$$e_{0r-i}(R(\bar I|\bar J)) = e_0(\bar I|\bar J) = e(\bar I,\bar A/0:\bar J^\infty) = e(I,A/(a_1,\ldots,a_i):J^\infty) > 0.$$
Since $e_i(I|J) = e_{0r-i}(R(\bar I|\bar J))$, this completes the proof of Theorem \ref{main}. \end{pf}

Theorem \ref{main} requires the existence of elements in $J$ with special properties. However, if the residue field of $A$ is infinite, we can always find such elements either by using Lemma \ref{exist2} or by the following notion.\sk 

Let $k = A/\frak m$ be the residue field of $A$ and $J = (c_1,\ldots,c_r)$. Following [Te] and [S] we say that a given property holds for a {\it sufficiently general element} $a \in J$ if there exists a non-empty Zariski-open subset $U \subseteq k^r$ such that whenever $a = \sum_{j=1}^r \alpha_jc_j$ and the image of $(\alpha_1,\ldots,\alpha_r)$ in $k^r$ belongs to $U$, then the given property holds for $a$.

\begin{Lemma} \label{general} Assume that $k$ is infinite. Let $a_1,\ldots,a_i$ be sufficiently general elements of $J$. Then their images in $J/IJ$ and $J/J^2$ form filter-regular sequences in $R(I|J)$ and $R(J|I)$. \end{Lemma}

\begin{pf}  Let $x_1,\ldots,x_i$ be the images of $a_1,\ldots,a_i$ in $J/IJ =  R_{(1,0)}$.  By definition, the filter-regularity of $x_1,\ldots,x_i$ means that $x_h \not\in \wp$ for every associated prime ideal $\wp \not\supseteq R_+$ of $(a_1,\ldots,a_{h-1})$, $h = 1,\ldots,i$. Since $R_+ \subset R_1$, $\wp \not\supseteq R_1$. Hence $\wp_{(1,0)} \neq R_{(1,0)} = J/IJ$. By Nakayama's Lemma $\wp_{(1,0)}$ maps to a  proper subspace $V(\wp)$ of the vector space $J/{\frak m}J$. Thus, if $a_h = \sum_{j=1}^r\alpha_jc_j$, there exists a non-empty Zariski-open subset $U \subseteq k^r$ such that whenever the images of $\alpha_1,\ldots,\alpha_r$ in $A/\frak m$ correspond to a point in $U$, then the image of $a_h$ in $J/{\frak m}J$ avoids $V(\wp)$, whence $x_h \not\in \wp$. Thus, $x_1,\ldots,x_i$ is a filter-regular sequence in $R(I|J)$. Similarly, we can show that the images of $a_1,\ldots,a_i$ in $J/J^2$ form a filter-regular sequence in $R(J|I)$. \end{pf}

As a consequence of Theorem \ref{main} we obtain the rigidity of mixed multiplicities (see Section 2 for the definition) and the independence of their positivity from the ideal $I$.

\begin{Corollary} \label{rigid3} Let $\rho = \max \{i|\ e_i(I|J) > 0\}.$ Then\par
{\rm (i) } $\height J -1 \le \rho < s(J)$,\par
{\rm (ii) } $e_i(I|J) > 0$ for $0 \le i \le \rho$,\par
{\rm (iii)}  $\max \{i|\ e_i(I'|J) > 0\} = \rho$ for any $\frak m$-primary ideal $I'$ of $A$. \end{Corollary}

\begin{pf} Without restriction we may assume that $k$ is infinite. By Lemma \ref{general} we can find elements $a_1,\ldots,a_{s(J)}$ in $J$ such that their images $x_1,\ldots,x_{s(J)}$ in $J/IJ$ and $x'_1,\ldots,x'_{s(J)}$ in $J/J^2$ form filter-regular sequences in $R(I|J)$ and $R(J|I)$, respectively. \par
(i) By Proposition \ref{deg u} we have $\rho < s(J)$. Let $s = \height J$. To prove $s - 1 \le \rho$ we may assume that $s \ge 1$. We can choose $a_1,\ldots,a_s$ such that $(a_1,\ldots,a_s)$ is part of a system of parameters of $A$. Then $J$ is not contained in any associated prime ideal $\frak p$ of $A$ with $\dim A/{\frak p} = \dim A$. From this it follows that such a  prime ideal is also an associated prime ideal of $0:J^\infty$. Hence $\dim A/0:J^\infty = \dim A$. Similarly, we can show that 
$$\dim A/(a_1,\ldots,a_{s-1}):J^\infty = \dim A/(a_1,\ldots,a_{s-1}) = \dim A-s+1 = \dim A/0:J^\infty-s+1.$$
By Theorem \ref{main} this implies $e_{s-1}(I|J)> 0$, whence $s-1 \le \rho$.\par
(ii) For a fixed $i \le \rho$ let $\bar A = A/(a_1,\ldots,a_i)$, $\bar I = I\bar A$ and $\bar J = J\bar A$. As in the proof of Theorem \ref{main} we can show that  $\deg P_{R(\bar I|\bar J)}(u,v) \le r-i$ and $e_\rho(I|J) = e_{\rho-ir-\rho}(R(\bar I|\bar J)),$ where $r = \dim A/0:J^\infty-2$. Since  $e_{\rho-ir-\rho}(R(\bar I|\bar J)) = e_\rho(I|J) > 0$, $\deg P_{R(\bar I|\bar J)}(u,v) = r-i$. Using Lemma \ref{char} and Lemma \ref{commute} we can see that $x_{i+1},\ldots,x_\rho$ and $x'_i,\ldots,x'_\rho$ form filter-regular sequences in $R(\bar I|\bar J)$ and $R(\bar J|\bar I)$. Now we can apply Theorem \ref{main} to $R(\bar I|\bar J)$ and obtain 
$$\dim A/(a_1,\ldots,a_\rho):J^\infty = \dim \bar A/0:\bar J^\infty = \dim A/(a_1,\ldots,a_i):J^\infty - \rho+i.$$
In particular, $\dim A/(a_1,\ldots,a_\rho):J^\infty = \dim A/0:J^\infty - \rho$. Hence 
$$\dim A/(a_1,\ldots,a_i):J^\infty = \dim A/0:J^\infty - i.$$ By Theorem \ref{main} this implies 
$e_i(I|J) > 0$. \par
(iii)  By Lemma \ref{general}, we may assume that the images of $a_1,\ldots,a_{s(J)}$ in $J/I'J$ and $J/J^2$ also form filter-regular sequences in $R(I'|J)$ and $R(J|I')$, respectively. By Theorem \ref{main}, the positivity of both $e_i(I|J)$ and $e_i(I'|J)$ depends only on $\dim A/(a_1,\ldots,a_i):J^\infty$. Hence $\rho = \max\{i|\ e_i(I'|J) > 0\}.$  \end{pf}

Katz and Verma claimed that $e_i(I|J) > 0$ for $1 \le i \le s(J)-1$ [KV, Theorem 2.7 (ii)]. Using Theorem \ref{main} we can easily construct counter-examples with $e_{s(J)-1}(I|J) = 0$.\ms

\noindent{\bf Example.} Let $A = k[[x_1,x_2,x_3,x_4]]/(x_1) \cap (x_2,x_3)$. Let $I$ be the maximal ideal of $A$ and $J = (x_1,x_4)A$. Then $F(J) \cong k[x_1,x_4]$. Hence $s(J) = \dim F(J) = 2$. Using Lemma \ref{char} we can verify that the images of $x_4$ in $J/IJ$ and $J/J^2$ are filter-regular elements in $R(I|J)$ and $R(J|I)$. We have $0:J^\infty = 0$ and $x_4A:J^\infty = (x_2,x_3,x_4)A$. Hence $\dim A/x_4A:J^\infty = 1 < \dim A/0:J^\infty - 1 = 2$. By Theorem \ref{main} this implies $e_1(I|J) = 0$. \ms

For the positivity of $e_{s(J)-1}(I|J)$ we obtain the following sufficient condition. 

\begin{Corollary} \label{rigid4} Assume that $A/0:J^\infty$ satisfies the first chain condition. Then $e_i(I|J)> 0$ for $0 \le i \le s(J)-1$. \end{Corollary}

\begin{pf}  By Corollary \ref{rigid3} we only need to show that $e_{s(J)-1}(R) > 0$. Without restriction we may assume that the residue field $k$ of $A$ is infinite. Using Lemma \ref{general} we can find elements $a_1,\ldots,a_{s(J)}$ in $J$ such that their images in $J/IJ$ and $J/J^2$ form filter-regular sequences in $R(I|J)$ and $R(J|I)$, respectively. Let $\frak p$ be a minimal associated prime ideal of $(a_1,\ldots,a_{s(J)-1}):J^\infty$ with 
$$\dim A/(a_1,\ldots,a_{s(J)-1}):J^\infty = \dim A/{\frak p}.$$ 
Since ${\frak p} \not\supseteq J$, $[(a_1,\ldots,a_{i-1}):a_i]_{\frak p} = (a_1,\ldots,a_{i-1})_{\frak p}$, $i = 1,\ldots,s(J)-1$. Thus, $a_1,\ldots,a_{s(J)-1}$ form a regular sequence in $A_{\frak p}$. This implies $\dim A_{\frak p} = s(J)-1$. Since $A$ satisfies the first chain condition and since $\height {\frak p}/0:J^\infty = \dim A_{\frak p}$,  
$$\dim A/{\frak p} = \dim A/0:J^\infty - \height {\frak p}/0:J^\infty = \dim A/0:J^\infty - s(J) + 1.$$
Therefore, applying Theorem \ref{main} we obtain $e_{s(J)-1}(R) > 0$. \end{pf} 

Another interesting consequence of Theorem \ref{main} is that for the computation of $e_i(I|J)$ we may replace $I$ and $J$ by their minimal reductions.

\begin{Corollary} \label{reduction} Let $I'$ and $J'$ be arbitrary reductions of $I$ and $J$, respectively. Then $e_i(I|J) = e_i(I'|J')$ for $i = 0,\ldots,r$. \end{Corollary}

\begin{pf} Without restriction we may assume that $k$ is infinite. Let $a_1,\ldots,a_i$ be sufficiently general elements of $J'$. Since $J'$ is a reduction of $J$,  $(J'+IJ)/IJ$ and $(J'+J^2)/J^2$ generate reductions of $R(I|J)_{(1)}$ and $R(J|I)_{(1)}$, respectively. Using this fact we can show, similarly as for Lemma \ref{general}, that  the images of $a_1,\ldots,a_i$ in $J/IJ$ and $J/J^2$ form filter-regular sequences in $R(I|J)$ and $R(J|I)$. Now, by Theorem \ref{main} we have
\begin{eqnarray*} e_i(I|J)  & = & e(I,A/(a_1,\ldots,a_i):J^\infty)\\ e_i(I'|J')& = & e(I',A/(a_1,\ldots,a_i):(J')^\infty). \end{eqnarray*} Since $J$ and $J'$ share the same radical,   $(a_1,\ldots,a_i):J^\infty = (a_1,\ldots,a_i):(J')^\infty$. Hence
$$e(I,A/(a_1,\ldots,a_i):J^\infty) = e(I,A/(a_1,\ldots,a_i):(J')^\infty) = e(I',A/(a_1,\ldots,a_i):(J')^\infty),$$ where the latter equality follows from the fact that $I'$ is a reduction of $I$ [NR]. \end{pf} 

\section{Applications}\sk

Let $(A,\frak m)$ be a local ring with infinite residue field and $d = \dim A$. Let $I$ be an $\frak m$-primary ideal and $J$ an arbitrary ideal of $A$. The main aim of this section is to show how to compute the mixed multiplicities $e_i(I|J)$. \sk

First we will derive from Theorem \ref{main} two formulae for $e_i(I|J)$, $i = 0,\ldots,\height J-1$. The first formula will be used later for explicit computations. The second formula generalizes Risler and Teissier's result on mixed multiplicities of $\frak m$-primary ideals [Te, Ch. 0, Proposition 2.1]. Moreover, it helps us interpret the Milnor numbers of general  linear sections of analytic hypersurfaces as mixed multiplicities (see the example below). 

\begin{Proposition} \label{parameter} Let $i = 0,\ldots,\height J-1$. Let $a_1,\ldots,a_i$ and $b_1,\ldots,b_{d-i}$ be sufficiently general elements in $J$ and $I$, respectively. Then
$$e_i(I|J) =  e(I,A/(a_1,\ldots,a_i)) = e((a_1,\ldots,a_i,b_1,\ldots,b_{d-i}),A).$$ \end{Proposition}

\begin{pf} By Theorem \ref{main} and Lemma \ref{general} we have
$$e_i(I|J) = e(I,A/(a_1,\ldots,a_i):J^\infty).$$
It is clear that $a_1,\ldots,a_i$ is part of a system of parameters of $A$. Let $\wp$ be an arbitrary associated prime of $(a_1,\ldots,a_i)$ with $\dim A/\wp = \dim A/(a_1,\ldots,a_i) = d-i$. Since $\height \wp \le i < \height J$, $\wp$ does not contain $J$. Therefore, $(a_1,\ldots,a_i)$ and $(a_1,\ldots,a_i):J^\infty$ share the same $\wp$-primary component. So we can conclude that
$$e(I,A/(a_1,\ldots,a_i):J^\infty) = e(I,A/(a_1,\ldots,a_i)).$$
It is also clear that $b_1,\ldots,b_{d-i}$ generate a reduction of $I/(a_1,\ldots,a_i)$. Therefore, 
$$e(I,A/(a_1,\ldots,a_i)) =  e((b_1,\ldots,b_{d-i}),A/(a_1,\ldots,a_i)).$$
Without restriction we may assume that $b_j \not\in {\frak p}$ for any associated prime ideal ${\frak p} \neq \frak m$ of $(a_1,\ldots,a_i,b_1,\ldots,b_{j-1})$, $j = 1,\ldots,d-i$. Then 
$$\displaylines{e((b_1,\ldots,b_{d-i}),A/(a_1,\ldots,a_i)) =   \ell(A/(a_1,\ldots,a_i,b_1,\ldots,b_{d-i})) -  \hspace{2cm} \cr \hspace{2cm} -\ell((a_1,\ldots,a_i,b_1,\ldots,b_{d-i-1}):b_{d-i}/(a_1,\ldots,a_i,b_1,\ldots,b_{d-i-1}))}$$
by [AB, Corollary 4.8]. On the other hand, we may also assume that $a_j \not\in \frak p$ for any associated prime ideal ${\frak p} \not\supseteq J$ of $(a_1,\ldots,a_{j-1})$, $j = 1,\ldots,i$. If $\frak p$ is an associated prime ideal of $(a_1,\ldots,a_{j-1})$ with $\dim A/{\frak p} \ge d-j$, then ${\frak p} \not\supseteq J$ because $d-j > d-\height J \ge \dim A/J$, whence $a_j \not\in \frak p$. So we can apply [AB, Corollary 4.8] and obtain
$$\displaylines{e((a_1,\ldots,a_i,b_1,\ldots,b_{d-i}),A) = \ell(A/(a_1,\ldots,a_i,b_1,\ldots,b_{d-i})) - \hspace{2cm} \cr \hspace{2cm} -\ell((a_1,\ldots,a_i,b_1,\ldots,b_{d-i-1}):b_{d-i}/(a_1,\ldots,a_i,b_1,\ldots,b_{d-i-1})).}$$
This implies
$$e((b_1,\ldots,b_{d-i}),A/(a_1,\ldots,a_i)) = e((a_1,\ldots,a_i,b_1,\ldots,b_{d-i}),A).$$
Summing up all equations we get the conclusion.  \end{pf}

\noindent{\bf Example.} {\it Analytic hypersurfaces.} Let $A = {\Bbb C}\{z_1,\ldots,z_n\}$ be the ring of convergent power series in $n$ variables. Let $f \in A$ be the equation of an analytic hypersurface $(X,x) \subset ({\Bbb C}^n,0)$. Let $J(f) = (\partial f/\partial z_1,\ldots,\partial f/\partial z_n)$ be the Jacobian ideal of $f$ in $A$. If $(X,x)$ is an isolated singularity, one calls $\mu(X,x) = \ell(A/J(f))$ the Milnor number of $(X,x)$. If $(X,x)$ is not an isolated singularity, let $s$ be the codimension of the singular locus of $X$. Let $H_i$ be a general $i$-plane in ${\Bbb C}^n$ passing through $x$, $i \le s$. For $i \le s$ Teissier [T1] proved that the Milnor number of $(X \cap H_i,x)$ is a constant and denoted it by $\mu^{(i)}(X,x)$. Let $a_1,\ldots,a_i$ and $b_1,\ldots,b_{n-i}$ be sufficiently general elements in $J(f)$ and $\frak m$, respectively. It is easily seen that 
\begin{eqnarray*} \mu^{(i)}(X,x) & = &\ell(A/(a_1,\ldots,a_i,b_1,\ldots,b_{n-i}))\\
& = & e((a_1,\ldots,a_i,b_1,\ldots,b_{n-i}),A). \end{eqnarray*}
By Proposition \ref{parameter} we have $e_i({\frak m}|J(f)) = e((a_1,\ldots,a_i,b_1,\ldots,b_{n-i}),A)$, $i = 0,\ldots,s-1$. Hence
$$\mu^{(i)}(X,x) = e_i({\frak m}|J(f)).$$
This fact was proved by Teissier only for analytic hypersurfaces with isolated singularities [Te, Ch. 0, Proposition 2.10].\ms

Next we will compute $e_s(I|J)$, $s = \height J$, in the case $s < s(J)$. Recall that $J$ is called {\it generically a complete intersection} if $\height {\frak p} = d-\dim A/J$ and $J_{\frak p}$ is generated by $\height {\frak p}$ elements for every associated prime ideal $\frak p$ of $J$ with $\dim A/{\frak p} = \dim A/J$.  

\begin{Proposition} \label{deviation} Let $J$ be an ideal of $A$ with $0 <\height J = s < s(J)$. Assume that $J$ is generically a complete intersection. Let $a_1,\ldots,a_s$ and $b_1,\ldots,b_{d-s}$ be sufficiently general elements in $J$ and $I$, respectively. Then
$$e_s(I|J) =  e(I,A/(a_1,\ldots,a_s)) - e(I,A/J).$$ \end{Proposition}

\begin{pf} From the condition that $J$ is generically a complete intersection we can deduce that $J_{\frak p} = (a_1,\ldots,a_s)_{\frak p}$ for every associated prime ideal $\frak p$ of $J$ with $\dim A/{\frak p} = \dim A/J$. Hence we can write 
$$(a_1,\ldots,a_s) = J \cap Q,$$
where $Q$ does not have any associated prime ideal $\frak p$ with $\dim A/{\frak p} = \dim A/J$ which is also an associated prime ideal of $J$. It follows that $(a_1,\ldots,a_s):J^\infty = Q:J^\infty$. If $\dim A/Q:J^\infty < d-s$, then $e_s(I|J) = 0$ by Theorem \ref{main} and $e(I,A/(a_1,\ldots,a_s)) = e(I,A/J)$ by the associativity formula of multiplicity. If $\dim A/Q:J^\infty = d-s$, then 
$$e_i(I|J) = e(I,A/Q:J^\infty)$$
by Theorem \ref{main}. Moreover, $Q$ must have at least one associated prime ideal $\wp$ with $\dim A/\wp = d-s$. Hence $\dim A/Q = d-s$. Since the associated prime ideals $\wp$ of $Q$ with  $\dim A/\wp = d-s$ do not contain $J$, the ideals $Q$ and $Q:J^\infty$ have the same minimal associated prime ideals $\wp$ with $\dim A/\wp = d-s$. Now, using the associativity formula of multiplicity we obtain
$$e(I,A/Q:J^\infty) = e(I,A/Q) = e(I,A/(a_1,\ldots,a_s)) - e(I,A/J).$$
The proof of Proposition \ref{deviation} is now complete. \end{pf}

For any local ring $(B,\frak n)$ let $e(B)$ denote the multiplicity of $B$ with respect to $\frak n$. By Theorem \ref{main} we have 
$e_i({\frak m}|J) = e(A/(a_1,\ldots,a_i):J^\infty),$ where $a_1,\ldots,a_i$ are sufficiently general elements in $J$. As we have seen in the above propositions, we may replace $e(A/(a_1,\ldots,a_i):J^\infty)$ by $e(A/(a_1,\ldots,a_i))$ if $i < \height J$ or by $e(A/(a_1,\ldots,a_i)) - e(A/J)$ if $i = \height J$ and $J$ is generically a complete intersection. The problem now  is to express $e(A/(a_1,\ldots,a_i))$ in terms of  better understood invariants of $A$ and $J$.\sk

The following result generalizes a result of Verma  in the case $J$ is an $\frak m$-primary ideal of a two-dimensional regular local ring [V1, Proof of Theorem 4.1].

\begin{Corollary} \label{e1} Let $(A,\frak m)$ be a regular local ring and $J$ an ideal with $\height J \ge 2$. Let $o(J)$ denote the order of $J$, that is, the largest integer $n$ such that $J \subseteq {\frak m}^n$. Then $$e_1({\frak m}|J) = o(J).$$ \end{Corollary}

\begin{pf} Let $a$ be a sufficiently general element in $J$. By Proposition \ref{parameter} we have $e_1({\frak m}|J) = e(A/(a)).$ But $e(A/(a)) = e(G/(a)^*)$, where $G = \oplus_{n\ge 0}{\frak m}^n/{\frak m}^{n+1}$ is the associated graded ring of $A$ and $(a)^*$ denotes the initial ideal of $(a)$ in $G$. Since $A$ is regular, $G$ is regular and $(a)^*$ is generated by the initial element $a^*$ of $a$ in $G$. Thus, $e(A/(a)) = \deg a^*.$ Since $a$ is a sufficiently general element of $J$, $\deg a^* = o(J)$. Hence $e_1({\frak m}|J) = o(J).$ \end{pf}

\noindent{\bf Example.} Let $A$ be the ring of convergent power series in $n$ variables. Let $f \in A$ be the equation of an analytic hypersurface $(X,x)$ in ${\Bbb C}^n$. We have seen in the preceding example that the Milnor number $\mu^1(X,x)$ is equal to $e_1({\frak m}|J(f))$.
By Corollary \ref{e1}, $e_1({\frak m}|J(f)) = o(J(f))$. It is easily seen from the definition of $o(J(f))$ that $o(J(f)) = m(X,x)-1$, where $m(X,x)$ denotes the multiplicity of $f$ at the origin. So we recover the formula $\mu^1(X,x) = m(X,x)-1$ of Teissier [Te, Ch. 0, 1.6 (2)]. \ms

If $A$ is a polynomial ring over a field and $J$ is a homogeneous ideal, we can often estimate $e_2({\frak m}|J)$, where $\frak m$ is now the maximal graded ideal of $A$. Note that the mixed multiplicities of homogeneous ideals can be defined similarly as in the local case.

\begin{Corollary} \label{e2} Let $A$ be a polynomial ring over a field and $J$ a homogeneous ideal with $\height J \ge 2$. Assume that $J$ contains two forms of the least possible degrees $c_1,c_2$ having no common factors. Then\par  
{\rm (i) } $e_2({\frak m}|J) = c_1c_2$ if $\height J \ge 3$,\par
{\rm (ii) } $e_2({\frak m}|J) = c_1c_2 - e(A/J)$ if $\height J = 2$ and $J$ is generically a complete intersection. \end{Corollary}

\begin{pf} We first compute $e(A/(a_1,a_2))$, where $a_1,a_2$ are sufficiently general elements in $J$. It is well-known that $e(A/(a_1,a_2)) = e(A/(a_1,a_2)^*)$, where $(a_1,a_2)^*$ denotes the initial ideal of $(a_1,a_2)$. Let $a_1^*, a_2^*$ be the initial elements of $a_1,a_2$. Without restriction we may assume that $\deg a_1^* = c_1, \deg a_2^* = c_2$ and that $a_1^*,a_2^*$ is a regular sequence in $A$. Then $(a_1,a_2)^* = (a_1^*,a_2^*)$. From this it follows that $e(A/(a_1,a_2)^*) = c_1c_2$. Now we apply Proposition \ref{parameter} and Proposition \ref{deviation} and obtain
$$e_2({\frak m}|J) = e(A/(x_1,x_2)) = c_1c_2$$
if $\height J \ge 3$, and
$$e_2({\frak m}|J) = e(A/(a_1,a_2))-e(A/J) = c_1c_2 - e(A/J)$$ 
if $\height J \ge 2$ and $J$ is generically a complete intersection. \end{pf}

The computation of $e_i({\frak m}|J)$ for $i \le \height J$ becomes easier if $J$ is a homogeneous ideal generated by forms of the same degree.

\begin{Corollary} \label{homogen} Let $A$ be a standard graded algebra over a field with maximal graded ideal $\frak m$. Let $J$ be a homogeneous ideal generated by forms of the same degree $c$ with $\height J = s$. Then
$$e_i({\frak m}|J) = c^ie(A)$$
for $i <s$. Moreover, if $J$ is generically a complete intersection with $s(J) \ge s+1$, then
$$e_s({\frak m}|J) = c^se(A) - e(A/J).$$ \end{Corollary}

\begin{pf} The statements follow from the Proposition \ref{parameter} and Proposition \ref{deviation} and from the fact that every sufficiently general element in $J$ is a homogeneous form of degree $c$. \end{pf}

The above corollaries, though simple, give us a complete picture on the mixed multiplicities $e_i({\frak m}|J)$  in many interesting cases. Before presenting some examples below we will sketch how one uses mixed multiplicities to compute the multiplicity of the Rees algebra and the degree of some embedding of rational varieties obtained by blowing up projective spaces. \sk

Let $R(J) = \oplus_{n\ge 0}J^nt^n$ be the Rees algebra of $J$ and $N$ the maximal graded ideal of $R(J)$. The computation of $e(R(J)_N)$ is usually hard. Until now, there have been few classes of ideals where one can compute the multiplicity of the Rees algebras explicitly, see e.g. [V1], [HTU], [Tr]. However,  Verma [V2, Theorem 1.4] showed that
$$e(R(J)_N) = \sum_{i=0}^s e_i({\frak m}|P).$$\sk

If $A = k[x_0,\ldots,x_n]$ is a polynomial ring in $n+1$ variables over a field $k$ and if $J$ is generated by forms of the same degree $c$, then $k[J_{c +1}]$ is isomorphic to the coordinate ring of a projective embedding of the rational $n$-fold obtained by blowing up ${\Bbb P}^n$ along the subscheme defined by $J$ [GGP, Theorem 2.1], [CH, Lemma 1.1]. If we consider $R(J)$ as a standard bigraded algebra with $R(J)_{(u,v)} = (J^ut^u)_{uc+v}$, then $R(J) \cong R({\frak m}|J)$ and $k[J_{c+1}]$ is the diagonal subalgebra $R(J)_\Delta = \oplus_{u \ge 0}R_{(u,u)}$ of $R(J)$ [STV]. The Hilbert polynomial of $R(J)_\Delta$ is given by the formula 
$$P_{R(J)_\Delta}(u) = P_R(u,u) = \sum_{i=0}^n{e_i({\frak m}|J) \over i!(n-i)!}u^n +\ \text{terms of degree $< n$}.$$
Hence we have the following relationship between the multiplicity of $R(J)_\Delta$ and the mixed multiplicities of $R$ [STV, Proposition 2.3]:
$$e(R(J)_\Delta) = \sum_{i =0}^n{n \choose i}e_i({\frak m}|J).$$ 
 
\noindent {\bf Example.}  {\it Ideals of sets of points in ${\Bbb P}^2$ with $h$-vectors of decreasing type}. In this case, $J = {\frak p_1} \cap \ldots \cap {\frak p_r} \subset A = k[x_0,x_1,x_2]$, where ${\frak p}_1,\ldots,{\frak p}_r$ are the defining prime ideals of the given points.
Recall that the $h$-vector of $J$ is the sequence $h_{A/J}(n)-h_{A/J}(n-1), \ n \ge 1,$ and that it is of decreasing type if it is strictly decreasing once it starts to decrease. One may view sets of points with $h$-vector of decreasing type as plane sections of curves in ${\Bbb P}^3$ [MR], [CO]. As a consequence, $J$ contains two forms of the least two possible degrees having no common divisor, say  $c_1 \le c_2$. Now we may apply Proposition \ref{parameter} ($i = 0$), Corollary \ref{e1} and Corollary \ref{e2} and obtain
$$e_0({\frak m}|J) = 1,\ e_1({\frak m}|J) = c_1,\ e_2({\frak m}|J) = c_1c_2 - r.$$
Since $s(J) \le \dim A = 3$, $e_i({\frak m}|J) = 0$ for $i \ge 3$. Therefore,
$$e(R(J)_N) = 1 + c_1 + c_1c_2 - r.$$
If $J$ is generated by forms of the same degree $c$, then
$$e(R(J)_\Delta) =  1 + 2c + c^2 - r.$$
Notice that if $J$ is the ideal of a set of ${c+1 \choose 2}$ points in ${\Bbb P}^2$ which do not lie on a curve of degree $d-1$, then $R(J)_\Delta$ is the coordinate ring of the Room surface in ${\Bbb P}^{3c+2}$ recently studied by Geramita and Gimigliano [GG]. \ms

\noindent {\bf Example.} {\it Defining prime ideals of curves in ${\Bbb P}^3$ which lie on the quadric $x_0x_4-x_1x_2$.} In this case, $J$ is a prime ideal in $A = k[x_0,x_1,x_2,x_3]$ with $J^n = J^{(n)}$ for all $n \ge 1$ [HH, Proposition 4.1]. By [HH, Theorem 2.5] this implies $s(J) \le \dim A - 1 = 3.$ By Proposition \ref{parameter}, Corollary \ref{e1} and Corollaries \ref{e2}, we have
$$e_0({\frak m}|J) = 1,\ e_1({\frak m}|J) = 2,\ e_2({\frak m}|J) = 2c - e(A/J),$$
where $c$ is the least degree of a form in $J/(x_0x_4-x_1x_2)$. Since $e_i({\frak m}|J) = 0$ for $i \ge 3$,
$$e(R(J)_N) = 3 + 2c - e(A/J).$$
If $J$ is generated only by quadrics, then
$$e(R(J)_\Delta) =  1 + 6 + 12 - 3e(A/J) = 19 - 3e(A/J).$$
A special case of this example is the monomial curve $(t_1^{b+c}:t_1^ct_2^b:t_1^bt_2^c:t_2^{b+c})$, where $b < c$, which was already dealt with by Raghavan and Verma [RV, 3.3]. \ms

\noindent {\bf Example.} {\it Homogeneous prime ideals of analytic deviation 1 which are generated by forms of the same degree in polynomial rings}. Recall that the analytic deviation of $J$ is the difference $s(J) - \height J$ [HH]. Let $J \subset A = k[x_0,\ldots,x_n]$, $s =\height J$ and let $c$ be the degree of the generators of $J$. By Proposition \ref{parameter} and Proposition \ref{deviation} we have
$$e_i({\frak m}|J) = c^i,\ i = 0,\ldots,s-1,\ e_s({\frak m}|J) = c^s - e(A/J).$$
Since $e_i({\frak m}|J) = 0$ for $i \ge s(J) = s+1$,  
\begin{eqnarray*}e(R(J)_N) & = & 1 +c + \cdots + c^s - e(A/J)\\
e(R(J)_\Delta) & = & 1 + nc + \cdots + {n \choose s}c^s - {n \choose s}e(A/J).\end{eqnarray*}
See [RV, 3.4] for the mixed multiplicities of another case of ideals of analytic deviation 1.\ms

\noindent{\bf Remark.} Raghavan and Verma [RV] used a technique of [HTU] which is similar to that of Gr\"obner bases to derive a formula for $e_i({\frak m}|J)$ when $J$ is generated by certain quadratic sequences (a generalization of $d$-sequences). Their formula expresses $e_i({\frak m}|J)$ in terms of the multiplicity $e(A/(a_1,\ldots,a_i):a_{i+1})$. Under their assumptions we always have $(a_1,\ldots,a_i):a_{j+1} = (a_1,\ldots,a_i):J^\infty$. Hence using Theorem \ref{main} we can also recover their formula. We leave the reader to check the details. \ms

Finally we will use Theorem \ref{main} to describe the degree of the St\"uckrad-Vogel cycles in intersection theory in terms of mixed multiplicities of ideals. Our result together with a recent result of Achilles and Manaresi [AM] which interprets the degree of the St\"uckrad-Vogel cycles as  mixed multiplicties of certain bigraded algebra (see below) provide a close relationship between mixed multiplicties and intersection theory. \sk

Let us first recall the definition of the St\"uckrad-Vogel cycles. Let $X$ and $Y$ be equidimensional subschemes of ${\Bbb P}_k^n$. Let $V$ be the ruled join variety of $X$ and $Y$ in ${\Bbb P}_{k(t)}^{2n+1} = \Proj k(t)[x_0,\ldots,x_n,y_0,\ldots,y_n]$, where $k(t) = k(t_{ij}|\ 1 \le i \le n+1, 0 \le j \le n)$ is a pure transcendental extension of $k$. Let $D$ be the diagonal subspace of ${\Bbb P}_{k(t)}^{2n+1}$ given by the equations $x_0-y_0 = \cdots = x_n-y_n = 0$. For $i = 1,\ldots,n+1$ let $h_i$ denote the divisor of $V$ given by the equation $\sum_{j=0}^nt_{ij}(x_j-y_j) = 0$.
Following the notations of [AM] we can define certain cycles $w_i$ and $v_i$ on $V$ inductively  as follows. First we put $w_0 = [V]$. If $w_{i-1}$ is defined for some $i \ge 1$, we decompose $w_{i-1} \cap h_i = v_i+w_i,$ where the support of $v_i$ lies in $D$ and $w_i$ has no components contained in $D$. It is clear that $\dim v_i = \dim V-i$. We call $v_1,\ldots,v_{n+1}$ the St\"uckrad-Vogel cycles of the intersection $X \cap Y$. These cycles were introduced in order to prove a refined Bezout's theorem for improper intersections [SV1], [Vo]. Van Gastel [Ga] showed later that their $k$-rational components correspond to the distinguished varieties of Fulton's intersection theory. \sk

Let $I_X$ and $I_Y$ denote the defining ideals of $X$ and $Y$ in $k[x_0,\ldots,x_n]$ and $k[y_0,\ldots,y_n]$, respectively. Then $A = k(t)[x_0,\ldots,x_n,y_0,\ldots,y_n]/(I_X,I_Y)$ is the  homogeneous coordinate ring of $V$. Let $\frak m$ be the maximal graded ideal of $A$ and $J = (x_0-y_0,\ldots,x_n-y_n)A$. Let $G_J(A)$ denote the associated graded ring of $A$ with respect to $J$. The associated graded ring $G_{\frak m}(G_J(A))$  of $G_J(A)$ with respect to the ideal ${\frak m}G_J(A)$ is a bigraded algebra over $k(u)$ with 
$$G_{\frak m}(G_J(A))_{(u,v)} = ({\frak m}^vJ^u+J^{u+1})/({\frak m}^{v+1}J^u+J^{u+1}).$$
Let $d = \dim A$. Achilles and Manaresi [AM, Theorem 4.1] proved that
$$\deg v_i = e_{d-i-1i-1}(G_{\frak m}(G_J(A))).$$
Note that they used the bidegree $(v,u)$ for $G_{\frak m}(G_J(A))$. \par
Using Theorem \ref{main} we can describe  $\deg v_i$ and therefore $e_{d-i-1i-1}(G_{\frak m}(G_J(A)))$ in terms of the simpler mixed multiplicities $e_i({\frak m}|J)$.

\begin{Theorem}\label{Vogel} With the above notations we have
$$\deg v_i = e_{i-1}({\frak m}|J) - e_i({\frak m}|J).$$ \end{Theorem}

\begin{pf} Let $a_i = \sum_{j=0}^nt_{ij}(x_j-y_j)$, $i = 1,\ldots,n+1$. We have
$$((a_1,\ldots,a_{i-1}):J^\infty,a_i):J^\infty = (a_1,\ldots,a_i):J^\infty.$$
This implies that $(a_1,\ldots,a_i):J^\infty$ is the intersection of the primary components of $((a_1,\ldots,a_{i-1}):J^\infty,a_i)$ whose associated prime ideals do not contain $J$. Hence we may write
$$((a_1,\ldots,a_{i-1}):J^\infty,a_i) = Q_i \cap ((a_1,\ldots,a_i):J^\infty),$$
where $Q_i$ is the intersection of the primary components of $((a_1,\ldots,a_{i-1}):J^\infty,a_i)$ whose associated prime ideals contain $J$. From the inductive definition of the cycles  $v_i$ we can deduce that
$$\deg v_i  =  e(A/Q_i).$$

Since $a_i$ is a generic linear combination of the generators of $J$ and since the associated prime ideals of $(a_1,\ldots,a_{i-1}):J^\infty$ can be defined over the field $k(t_{hj}|\ 1 \le h \le i-1,\ 0 \le j \le n)$ and do not contain $J$, $a_i$ is a non-zerodivisor in $A/(a_1,\ldots,a_{i-1}):J^\infty$. Therefore
\begin{eqnarray*} e(A/(a_1,\ldots,a_{i-1}):J^\infty) & = & e(A/((a_1,\ldots,a_{i-1}):J^\infty,a_i))\\ & = & e(A/Q_i) + e(A/(a_1,\ldots,a_i):J^\infty). \end{eqnarray*}
>From this it follows that
$$\deg v_i = e(A/(a_1,\ldots,a_{i-1}):J^\infty)-e(A/((a_1,\ldots,a_{i-1}):J^\infty,a_i)).$$
Since the images of $a_1,\ldots,a_{n+1}$ in $J/{\frak m}J$ resp.~$J/J^2$ are generic linear combinations of the generators of $J/{\frak m}J$ resp.~$J/J^2$, they form filter-regular sequences in $R(I|J)$ resp.~$R(J|I)$. Therefore we may apply Theorem \ref{main} and obtain
\begin{eqnarray*} e_{i-1}({\frak m}|J) & = & e(A/(a_1,\ldots,a_{i-1}):J^\infty),\\
e_i({\frak m}|J) & = & e(A/(a_1,\ldots,a_i):J^\infty). \end{eqnarray*}
Putting these relations into the above formula for $\deg v_i$ we get\par 
$\hspace{5cm} \deg v_i = e_{i-1}({\frak m}|J) - e_i({\frak m}|J).$ \end{pf}\sk

\section*{References}\sk

\noindent [AM] R. Achilles and M. Manaresi, Multiplicities of a bigraded ring and intersection theory, Math. Ann. 309 (1997), 573-591. \par
\noindent [AB] M. Auslander and D.A. Buchsbaum, Codimension and multiplicity, Ann. Math. 68 (1958), 625-657.\par
\noindent [B] P.B. Bhattacharya, The Hilbert function of two ideals, Proc. Cambridge Phil. Soc. 53 (1957), 568-575.\par
\noindent [BR] M. Brodmann and J. Rung, Local cohomology and the connectedness dimension in algebraic varieties, Comment. Math. Helvetici 61 (1986), 481-490. \par
\noindent [BS] M. Brodmann and R.Y. Sharp, Local cohomology, an algebraic introduction with geometric applications, Cambridge University Press, 1998. \par
\noindent [CO] L. Chiantini and F. Orrechia, Plane sections of arithmetically normal curves in ${\Bbb P}^3$, in: Algebraic curves and projective geometry (Trento, 1988), 32-42, Lect. Notes in Math. 1389, Springer, 1989.\par
\noindent [CH] S. D. Cutkosky and J. Herzog, Cohen-Macaulay coordinate rings of blowup schemes, Comment. Math. Helvetici 72 (1997), 605-617. \par
\noindent [E] D. Eisenbud, Commutative algebra with a view toward algebraic geometry,  Springer-Verlag, 1995.\par
\noindent [Ga] L.J. van Gastel, Excess intersection and a correpondence principle, Invent. Math. 103 (1991), 197-221.\par
\noindent [GG] A.V. Geramita and A. Gimigliano, Generators for the defining ideal of certain rational surfaces, Duke Math. J. 62 (1991), 61-83. \par
\noindent [GGP] A.V. Geramita, A. Gimigliano, Y. Pitteloud, Graded Betti numbers of some embedded rational $n$-folds, Math. Ann. 301 (1995), 363-380. \par
\noindent [HHRT] M. Herrmann, E. Hyry, J. Ribbe and Z.M. Tang, Reduction numbers and multiplicities of multigraded structures, J. Algebra 197, 311-341 (1997). \par
\noindent [HTU] J. Herzog, N.V. Trung and B. Ulrich,  On the multiplicity of Rees algebras and associated graded rings of d-sequences,  J. Pure Appl. Algebra~80 (1992), 273-297. \par
\noindent [HH] S. Huckaba and C. Huneke, Powers of ideals having small analytic deviation, Amer. J. Math. 114 (1992), 367-403.\par
\noindent [KMV] D. Katz, S. Mandal and J. Verma, Hilbert function of bigraded algebras, in:  A. Simis, N.V. Trung and G. Valla  (eds.), Commutative Algebra (ICTP, Trieste, 1992), 291-302, World Scientific, 1994.\par
\noindent [KV] D. Katz and J.K. Verma, Extended Rees algebras and mixed multiplicities, Math. Z. 202 (1989), 111-128.\par
\noindent [KR1] D. Kirby and D. Rees, Multiplicities in graded rings I: the general theory, Contemporary Mathematics 159 (1994), 209-267. \par
\noindent [KR2] D. Kirby and D. Rees, Multiplicities in graded rings II: integral equivalence and the Buchsbaum-Rim multiplicity, Math. Proc. Cambridge Phil. Soc. 119 (1996), 425-445.\par
\noindent [KT1] S. Kleiman and A. Thorup, A geometric theory of the Buchsbaum-Rim multiplicity, J. Algebra 167 (1994), 168-231.\par
\noindent [KT2] S. Kleiman and A. Thorup, Mixed Buchsbaum-Rim multiplicities, Amer. J. Math. 118 (1996), 529-569.\par
\noindent [MR] R. Maggioni and A. Ragusa, The Hilbert functions of generic plane sections of curves in ${\Bbb P}^3$, Invent. Math. 91 (1988), 253-258. \par
\noindent [NR] D.G. Northcott and D. Rees, Reductions of ideals in local rings, Proc. Cambridge Phil. Soc. 50 (1954), 145-158. \par
\noindent [RV] K.N. Raghavan and J.K. Verma, Mixed Hilbert coefficients of homogeneous $d$-sequences and quadratic sequences, J. Algebra 195 (1997), 211-232. \par
\noindent [R1] D. Rees, $\frak a$-transforms of local rings and a theorem on multiplicities of ideals, Proc. Cambridge Philos. Soc. 57 (1961), 8-17. \par
\noindent [R2] D. Rees, Generalizations of reductions and mixed multiplicities, J. London Math. Soc. 29 (1984), 423-432.\par
\noindent [Ro] P. Roberts, Local Chern classes, multiplicities and perfect complexes, M\'emoire Soc. Math. France 38 (1989), 145-161. \par
\noindent [STV] A. Simis, N.V. Trung and G. Valla, The diagonal subalgebras of a blow-up ring, J. Pure Appl. Algebra 125 (1998), 305-328. \par
\noindent [SV1] J. St\"uckrad and W. Vogel, An algebraic approach to the intersection theory, in: The curves seminar Vol. II, 1-32, Queen's papers in pure and applied mathematics, No. 61, Kingston, 1982. \par
\noindent [SV2] J. St\"uckrad and W. Vogel, Buchsbaum rings and applications, VEB Deutscher Verlag der Wisssenschaften, Berlin, 1986.\par
\noindent [Sw] I. Swanson, Mixed multiplicities, joint reductions and quasi-unmixed local rings, J. London Math. Soc. 48 (1993), no. 1, 1-14.\par
\noindent [Te] B. Teissier, Cycles \'evanescents, sections planes, et conditions de Whitney, Singularit\'es \`a Carg\`ese 1972, Ast\`erisque 7-8 (1973), 285-362. \par
\noindent [Tr] N.V. Trung,  Filter-regular sequences and multiplicity of blow-up rings of ideals of the principal class,  J. Math. Kyoto Univ.~33 (1993), 665-683.\par
\noindent [V1] J.K. Verma, Rees algebras and mixed multiplicities, Proc. Amer. Math. Soc. 104 (1988), 1036-1044. \par
\noindent [V2] J.K. Verma, Multigraded Rees algebras and mixed multiplicities, J. Pure Appl. Algebra 77 (1992), 219-228. \par
\noindent [Vo] W. Vogel, Lectures on Bezout's theorem, Tata Institute of Fundamental Research Lect. Notes 74, Springer, 1984.\par
\noindent [W] B.L. van der Waerden, On Hilbert's function, series of composition of ideals and a generalization of the theorem of Bezout, Proc. K. Akad. Wet. Amsterdam 31 (1928), 749-770.\par
\end{document}